# A curious instability phenomenon for a rounded corner in presence of a negative material

Lucas Chesnel[1], Xavier Claeys[2], Sergey A. Nazarov[3]

[1] Department of Mathematics and Systems Analysis, Aalto University, P.O. Box 11100, FI-00076 Aalto, Finland;
[2] Laboratory Jacques Louis Lions, University Pierre et Marie Curie, 4 place Jussieu, 75005 Paris, France;
[3] Faculty of Mathematics and Mechanics, St. Petersburg State University, Universitetsky prospekt, 28, 198504, Peterhof, St. Petersburg, Russia;

E-mail: lucas.chesnel@aalto.fi, claeys@ann.jussieu.fr, srgnazarov@yahoo.co.uk

**Abstract.** We study a 2D scalar harmonic wave transmission problem between a classical dielectric and a medium with a real-valued negative permittivity/permeability which models a metal at optical frequency or an ideal negative metamaterial. We highlight an unusual instability phenomenon for this problem when the interface between the two media presents a rounded corner. To establish this result, we provide an asymptotic expansion of the solution, when it is well-defined, in the geometry with a rounded corner. Then, we prove error estimates. Finally, a careful study of the asymptotic expansion allows us to conclude that the solution, when it is well-defined, depends critically on the value of the rounding parameter. We end the paper with a numerical illustration of this instability phenomenon.

## 1 Introduction

In electromagnetics, recent years have seen a growing interest in the use of negative materials in the development of new technological devices. These negative materials owe their name to the fact that they can be modeled for certain ranges of frequencies, neglecting the dissipation, by real negative physical parameters (permittivity $\varepsilon$ and/or permeability $\mu$). To summarize, they are divided into two major families. The *negative metamaterials*, and in particular the left-handed materials for which we have both $\varepsilon < 0$ and $\mu < 0$, are complex structures made of small resonators, chosen so that the macroscopic medium behave as if its physical parameters were negative. For a mathematical justification of the homogenization process, we refer the reader for example to [12, 13, 14]. *Metals* in visible range constitute the second family of negative materials. They are used especially in plasmonic technologies [2, 16, 53, 22] which would allow important advances in miniaturization. In this context, a key issue is to be able to manipulate light and in particular, to focus energy in specific areas of space. To do this, physicists use metallic devices with corners and edges [49, 3, 41].

Because of the sign-changing of the physical parameters, the study of time harmonic Maxwell's equations in devices involving negative materials raises challenging questions both from a theoretical and a numerical perspective [45, 43, 21]. Using a variational approach, it has been proved in [8, 10, 4] that the scalar problem equivalent to Maxwell's equations in 2D configurations, turns out to be of Fredholm type in the classical functional framework whenever the contrast (ratios of the values of $\varepsilon$ or $\mu$ across the interface) lies outside some interval, which always contains the value $-1$. Moreover, this interval reduces to $\{-1\}$ if and if only the interface between the positive material and the negative material is smooth (of class $\mathscr{C}^1$). Analogous results have been obtained by techniques of boundary integral equations in [19]. The numerical approximation of the solution of this scalar problem, based on classical finite element methods, has been investigated in [10, 40, 18]. Under some assumptions on the meshes, the discretized problem is well-posed and its solution converges to the solution of the continuous problem. The study of Maxwell's equations has been carried out in [9, 6, 5, 17]. The influence of corners of the interface, studied in [33, 50, 52, 25, 44, 23], has been clarified in [7] for the scalar problem (see also the previous works [20, 11, 46] where the general theory [26, 32, 35, 27] is extended to this configuration where the operator is not strongly elliptic). In [7], following [35, 38, 1, 39], the authors prove that when the contrast of the physical parameters lies inside the critical interval, Fredholm property is lost because of the existence of two strongly oscillating singularities at the corner. In such a case, Fredholmness can be recovered by adding to the functional framework one of the two singularities, selected by means of a limiting absorption principle, and by working in a special weighted Sobolev setting with weight centered at the corner. This functional framework amounts to prescribing a radiation condition at the corner.

Such a special functional framework seems an uncomfortable situation though, at least from a physical point of view. Indeed, it leads to work with solutions which are not of finite energy (their $H^1$-norm is infinite). We can imagine at least two ways to regularize this problem. We could for example add some dissipation to the medium. It has been proven in [7] (this is the limiting absorption principle mentioned above) that the solution of the problem with a small dissipation converges to the solution of the limit problem, of course in a weaker norm than the $H^1$-norm. Another possible regularization that may appear natural would consist in considering slightly rounded corners, instead of real corners at the interface. In the present article, we examine the validity of such a regularization process, studying the convergence of our diffusion problem in the case where the interface contains a rounded corner that tends (in the geometrical sense) to a sharp corner. Quite unexpectedly, we prove an unusual instability phenomenon: the solution, in the configuration where the corner of the interface is slightly rounded, depends critically on the value of the rounding parameter and does not converge when this parameter goes to zero.

The outline of the paper is the following. In Section 2, we describe in detail the problem and the geometry that we want to consider, namely a diffusion equation with a sign-changing coefficient in a cavity with an interface containing a rounded corner close to the boundary. This geometry depends on a small parameter $\delta$ corresponding to the rounding parameter of the corner. In Section 3, we consider two problems set in limit geometries related to the initial problem. Using results of [7], we state important properties of such problems and describe their solutions in terms of asymptotic expansions. In Section 4, we propose an asymptotic expansion with respect to $\delta$ for the solution to the problem that we have described in Section 2. This expansion is built using matched asymptotics, in accordance with the usual procedure [29], [31, Chap. 4,5]. In Section 5, we show that the norm of the resolvent of the problem considered in Section 2 admits at most a logarithmic growth with respect to $\delta$. The construction of an "asymptotic inverse" which appears in the approach is again due to [29], [31, Chap. 4] (the reader might also find a simple presentation of this method for concrete problems in mathematical physics in [36, 37, 15]). In Section 6, we state and prove the important Theorem 6.1, which provides an error estimate between the solution to the initial problem and the first two terms of its asymptotic expansion with respect to $\delta$. Thanks to this theorem, we show the main result of the paper, namely the instability phenomenon with respect to the rounded corner. We conclude with numerical experiments illustrating this instability phenomenon.

## 2 Description of the problem

We consider a domain $\Omega \subset \mathbb{R}^2$, *i.e.* a bounded and connected open subset of $\mathbb{R}^2$ with Lipschitz boundary $\partial\Omega$ (see Fig. 1 below). We assume that $\Omega$ is partitioned into two sub-domains $\Omega_\pm^\delta$ so that $\overline{\Omega} = \overline{\Omega_+^\delta} \cup \overline{\Omega_-^\delta}$ with $\Omega_+^\delta \cap \Omega_-^\delta = \emptyset$. We consider a straight segment $\Sigma^0$ that intersects $\partial\Omega$ at only two points $O$ and $O'$. We assume that $\partial\Omega$ is straight in a neighbourhood of $O$, $O'$, and that at $O'$, $\Sigma^0$ is perpendicular to $\partial\Omega$. We also assume that the interface $\Sigma^\delta := \overline{\Omega_+^\delta} \cap \overline{\Omega_-^\delta}$ coincides with $\Sigma^0$ outside the disk $\mathrm{D}(O,\delta)$.

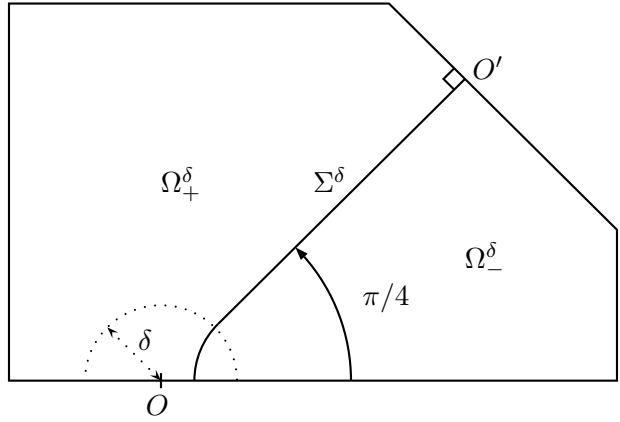

Fig. 1: Geometry of the problem.

In the sequel, we shall denote by $(r,\theta)$ the polar coordinates centered at $O$ such that $\theta = 0$ or $\pi$ at the boundary in a neighbourhood of $O$. As $\delta \to 0$, the sub-domains $\Omega_\pm^\delta$ turn into $\Omega_\pm^0$ and we assume that there exists a disk $\mathrm{D}(O, r_0)$ centered at $O$ such that $\Omega_-^0 \cap \mathrm{D}(O, r_0) = \{(r\cos\theta, r\sin\theta) \in \mathbb{R}^2 \mid 0 < r < r_0,\ 0 < \theta < \pi/4\}$ and $\Omega_+^0 \cap \mathrm{D}(O, r_0) = \{(r\cos\theta, r\sin\theta) \in \mathbb{R}^2 \mid 0 < r < r_0,\ \pi/4 < \theta < \pi\}$. We consider the value $\pi/4$ for the opening of the corner simply because it allows explicit calculus (see §3.1). There is no difficulty to adapt the rest of the forthcoming analysis for other values of this angle. To fix our ideas, and without restriction, we assume that we can take $r_0 = 2$, i.e. that $(\mathrm{D}(O,2) \cap \mathbb{R} \times \mathbb{R}_+^*) \subset \Omega$ with $\mathbb{R} \times \mathbb{R}_+^* = \{(x,y) \in \mathbb{R}^2 \mid y > 0\}$.

## 2.1 Geometry of the rounded corner

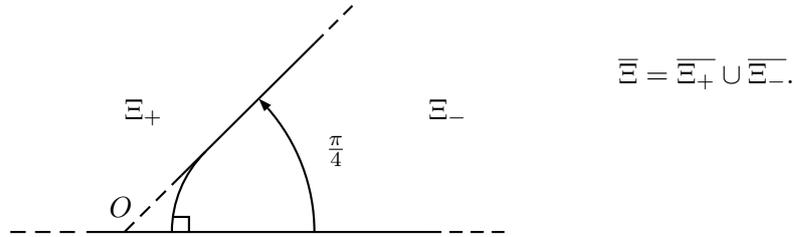

$$\overline{\Xi} = \overline{\Xi_+} \cup \overline{\Xi_-}.$$

Fig. 2: Frozen geometry.

The set $\Sigma^\delta \cap \mathrm{D}(O,\delta)$ will be defined as follows. Let $\Xi := \mathbb{R} \times \mathbb{R}_+^*$ refer to the upper half plane partitioned by means of two open sets $\Xi_\pm$ such that $\overline{\Xi} = \overline{\Xi_+} \cup \overline{\Xi_-}$ and $\Xi_+ \cap \Xi_- = \emptyset$. We assume that $\Gamma := \overline{\Xi_+} \cap \overline{\Xi_-}$ is a curve $\Gamma = \{\varphi_\Gamma(t),\ t \in [0, +\infty)\}$ where $\varphi_\Gamma$ is a $\mathscr{C}^\infty$ function such that $\partial_t \varphi_\Gamma(0)$ is orthogonal to the $x$-axis and $\varphi_\Gamma(t) = (t,t)$ for $t \geq 1$, see Fig. 2 below. In the neighbourhood of the corner, we assume that $\Omega_\pm^\delta$ can be defined from $\Xi_\pm$ by self similarity:

$$\Omega_\pm^\delta \cap \mathrm{D}(O,\delta) = \{\ \mathbf{x} \in \mathbb{R}^2 \mid \mathbf{x}/\delta \in \Xi_\pm \cap \mathrm{D}(O,1)\ \}.$$

**Remark 2.1.** *We use the word "corner" to refer to the geometry of the interface. It would have been more precise to describe it as a "half corner" (see Fig. 3). But the important point is that this problem shares the same specificities as "real corner" problems. We focus on this configuration only to simplify the presentation.*

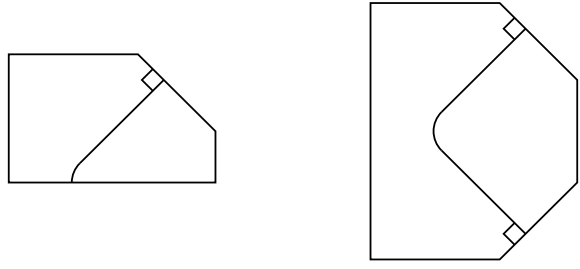

Fig. 3: "Half rounded corner" on left, "real rounded corner" on right.

## 2.2 The problem under study

First of all, let us fix some basic notations. In the sequel, for any open subset $\omega \subset \mathbb{R}^d$ with $d = 1, 2$, the space $\mathrm{L}^2(\omega)$ will refer to the Lebesgue space of square integrable functions equipped with the scalar product $u, v \mapsto (u, v)_\omega$, where we set

$$(u, v)_\omega = \int_\omega u\, \overline{v}\, d\mathbf{x}, \quad \forall u, v \in \mathrm{L}^2(\omega).$$

We denote $\|v\|_\omega := \sqrt{(v, v)_\omega}$. We will consider the Sobolev space $\mathrm{H}^1(\omega) := \{v \in \mathrm{L}^2(\omega) \mid \nabla v \in \mathrm{L}^2(\omega)\}$, and define $\mathrm{H}^1_0(\omega) := \{v \in \mathrm{H}^1(\omega) \mid v|_{\partial \omega} = 0\}$ equipped with the norm $\|v\|_{\mathrm{H}^1_0(\omega)} := \|\nabla v\|_\omega$. The space $\mathrm{H}^{-1}(\omega)$ will refer to the topological dual to $\mathrm{H}^1_0(\omega)$, made of the antilinear forms on $\mathrm{H}^1_0(\omega)$. It will be endowed with the intrinsic norm

$$\|f\|_{\mathrm{H}^{-1}(\omega)} := \sup_{v \in \mathrm{H}^1_0(\omega) \setminus \{0\}} \frac{|\langle f, v \rangle_\omega|}{\|v\|_{\mathrm{H}^1_0(\omega)}}, \quad \forall f \in \mathrm{H}^{-1}(\omega),$$

where $\langle \,,\, \rangle_\omega$ refers to the duality pairing between $\mathrm{H}^{-1}(\omega)$ and $\mathrm{H}^1_0(\omega)$.

The present article will focus on a transmission problem with a sign-changing coefficient. Define the function $\sigma^\delta : \Omega \to \mathbb{R}$ by $\sigma^\delta = \sigma_\pm$ in $\Omega^\delta_\pm$, where $\sigma_+ > 0$ and $\sigma_- < 0$ are constants. For $f \in \mathrm{H}^{-1}(\Omega)$, consider the problem, in the sense of distributions:

$$\left| \begin{array}{l} \text{Find } u^\delta \in \mathrm{H}^1_0(\Omega) \text{ such that} \\ -\mathrm{div}(\sigma^\delta \nabla u^\delta) = f \quad \text{in } \Omega. \end{array} \right. \quad (1)$$

Notice that when $f \in \mathrm{L}^2(\Omega)$, this problem also writes

$$\left| \begin{array}{rcll} \text{Find } (u^\delta_+, u^\delta_-) \in \mathrm{H}^1(\Omega^\delta_+) \times \mathrm{H}^1(\Omega^\delta_-) \text{ such that} \\ -\sigma_\pm \Delta u^\delta_\pm &=& f_\pm & \text{in } \Omega^\delta_\pm \\ u^\delta_+ - u^\delta_- &=& 0 & \text{on } \Sigma^\delta \\ \sigma_+ \partial_{\nu^\delta} u^\delta_+ - \sigma_- \partial_{\nu^\delta} u^\delta_- &=& 0 & \text{on } \Sigma^\delta \\ u^\delta_\pm &=& 0 & \text{on } \partial\Omega^\delta_\pm \cap \partial\Omega, \end{array} \right.$$

where $\nu^\delta$ denotes the unit outward normal vector to $\Sigma^\delta$ orientated from $\Omega^\delta_+$ to $\Omega^\delta_-$. Problem (1) can be reformulated as the integral identity

$$\left| \begin{array}{l} \text{Find } u^\delta \in \mathrm{H}^1_0(\Omega) \text{ such that} \\ (\sigma^\delta \nabla u^\delta, \nabla v)_\Omega = \langle f, v \rangle_\Omega, \quad \forall v \in \mathrm{H}^1_0(\Omega). \end{array} \right. \quad (2)$$

This leads us to introduce the continuous linear operator $\mathrm{A}^\delta : \mathrm{H}^1_0(\Omega) \to \mathrm{H}^{-1}(\Omega)$ defined by

$$\langle \mathrm{A}^\delta u, v \rangle_\Omega = (\sigma^\delta \nabla u, \nabla v)_\Omega, \quad \forall u, v \in \mathrm{H}^1_0(\Omega).$$

The function $u^\delta$ satisfies Problem (1) if and only if it satisfies $\mathrm{A}^\delta u^\delta = f$. As it is known from [4], $\mathrm{A}^\delta$ is a Fredholm operator of index 0[1] (with possibly a non trivial kernel) whenever $\kappa_\sigma := \sigma_-/\sigma_+ \neq -1$, as the interface $\Sigma^\delta$ is smooth and meets $\partial\Omega$ orthogonally. In [4, Th. 6.2], it is also proved that, as soon as $\Sigma^\delta$ presents a straight section, in the case $\kappa_\sigma = \sigma_-/\sigma_+ = -1$, the operator $\mathrm{A}^\delta$ is not of Fredholm type. Actually, for this configuration, one can check that ellipticity is lost for Problem (1) (see [48, 47] and [28]). Therefore, the situation $\kappa_\sigma = -1$ cannot be studied with the tools we propose. We refer the reader to [42] for more details concerning this case and we discard it from now on.

In the present paper, we are interested in studying the behaviour, as $\delta \to 0$, of the solution $u^\delta$ to the Problem (1) when it is well-defined. Note that for $\delta = 0$, the interface no longer meets $\partial\Omega$ perpendicularly.

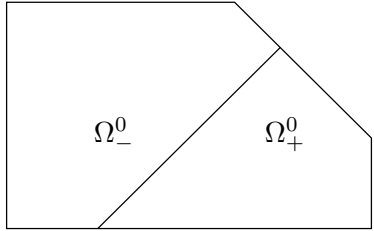

Fig. 4: Geometry for $\delta = 0$.

As shown in [4] and as mentioned in the introduction, there exist values of the contrasts $\kappa_\sigma = \sigma_-/\sigma_+$ for which the operator $\mathrm{A}^0$ fails to be of Fredholm type, because of the existence of two strongly oscillating singularities at the corner point $O$. More precisely, for the present geometrical configuration, $\mathrm{A}^0$ is a Fredholm operator if and only if, $\kappa_\sigma \in \mathbb{R}_-^* := (-\infty, 0)$ satisfies $\kappa_\sigma \notin [-1, -1/3]$. Here, the value 3 comes from the ratio of the two apertures: $3 = (\pi - \pi/4)/(\pi/4)$.

When $\mathrm{A}^0$ is of Fredholm type, there is no qualitative difference between Problem (1) for $\delta > 0$, and Problem (1) for $\delta = 0$. In this case, using the analysis we provide in this article (and which was introduced in [29], [31, Chap. 4]) we can prove that, if $\mathrm{A}^0$ is injective, then $\mathrm{A}^\delta$ is injective for $\delta$ small enough. Moreover, defining $u^\delta = (\mathrm{A}^\delta)^{-1}f$ and $u = (\mathrm{A}^0)^{-1}f$, we can show that the sequence $(u^\delta)$ converges to $u$ for the norm $\|\ \|_{\mathrm{H}^1(\Omega)}$. Since these results can be obtained from the approach we develop here, in a more classical way, we have chosen not to present them.

When $\mathrm{A}^0$ is not of Fredholm type, there is a qualitative difference between Problem (1) for $\delta > 0$, and Problem (1) for $\delta = 0$. The purpose of the present document is to study such a qualitative transition. When $\kappa_\sigma = -1/3$, the singularities associated to the corner have a more complex structure, with a logarithmic term. The approach we present in this paper needs to be adapted to deal with this configuration (see [29], [31, Chap. 4]). In the sequel, we discard this limit case, and therefore, we assume that

$$\kappa_\sigma = \sigma_-/\sigma_+ \in (-1, -1/3). \qquad (3)$$

The important message here is that under Condition (3), for all $\delta > 0$, $\mathrm{A}^\delta$ is of Fredholm type but $\mathrm{A}^0$ is not.

---

[1] All through the paper, each time the term "Fredholmness" is used for $\mathrm{A}^\delta$ or $\mathrm{A}^0$, it is understood that these operators act from $\mathrm{H}_0^1(\Omega)$ to $\mathrm{H}^{-1}(\Omega)$.

## 3 Limit geometries

As is usual in asymptotic analysis, we consider problems set in limit geometries independent from $\delta$ to study Problem (1) as $\delta \to 0$. The present analysis needs to consider a far field problem set in an outer limit geometry, and a near field problem set in an inner limit geometry.

### 3.1 Outer limit geometry

The first limit geometry that we wish to consider is obtained from Fig. 1 for $\delta = 0$, see Fig. 4 above. We call it the *outer limit geometry*. Since in this geometry, the interface $\Sigma^0 = \overline{\Omega^0_-} \cap \overline{\Omega^0_+}$ does not meet $\partial\Omega$ perpendicularly, the operator $A^0$ is not of Fredholm type when $\kappa_\sigma \in (-1, -1/3)$. The goal of this paragraph is to define an appropriate functional framework for this problem recalling the results obtained in [7]. For some $g$, we wish to consider a problem of the form

$$\text{Find } v \text{ such that } \quad -\text{div}(\sigma^0 \nabla v) = g \quad \text{in } \Omega, \quad v = 0 \text{ on } \partial\Omega. \tag{4}$$

Here, $\sigma^0$ is the function such that $\sigma^0 = \sigma_\pm$ in $\Omega^0_\pm$. We did not specify any functional framework: we have to introduce adapted weighted spaces that endow the problem above with a Fredholm structure.

In order to study Problem (4), according to the Kondratiev theory, we need first to describe the singularities associated to the corner point $O$. Once singularities at $O$ have been computed, all the results become a consequence of the general theory of [26, 32] (see also [35, 27]). Singularities are functions of separate variables in polar coordinates which satisfy the homogeneous problem in the infinite corner. According to §4.1 in [7], the problem of finding couples $(\lambda, \varphi) \in \mathbb{C} \times H^1_0(0, \pi)$ such that $\text{div}(\sigma^0 \nabla(r^\lambda \varphi(\theta))) = 0$ in $\Omega$ has non-trivial solutions only for $\lambda$ belonging to the following set of singular exponents

$$\Lambda := (2\mathbb{Z} \setminus \{0\}) \cup \{i\mu + 4\mathbb{Z}\} \cup \{-i\mu + 4\mathbb{Z}\}$$

$$\text{with} \quad \mu = -\frac{2}{\pi} \ln \left[ \frac{1}{2} \frac{\sigma_- - \sigma_+}{\sigma_- + \sigma_+} + i\sqrt{1 - \left(\frac{1}{2} \frac{\sigma_- - \sigma_+}{\sigma_- + \sigma_+}\right)^2} \right]. \tag{5}$$

In the case where $\sigma_-/\sigma_+ \in (-1, -1/3)$, we have $\mu \in (0, +\infty)$, so that the set $\Lambda$ contains only two points in the strip $|\Re e\{\lambda\}| < 2$, namely $\pm i\mu$. For $\lambda = \pm i\mu$, the space of functions $\varphi \in H^1_0(0, \pi)$ such that $\text{div}(\sigma^0 \nabla(r^{\pm i\mu} \varphi(\theta))) = 0$ is one dimensional, and generated by some $\phi(\theta)$ (both for $i\mu$ and $-i\mu$), see [7, §4.1],

$$\phi(\theta) = c_\phi \frac{\sinh(\mu\theta)}{\sinh(\mu\pi/4)} \quad \text{on } [0, \pi/4], \quad \text{and} \quad \phi(\theta) = c_\phi \frac{\sinh(\mu(\pi-\theta))}{\sinh(\mu \, 3\pi/4)} \quad \text{on } [\pi/4, \pi] \tag{6}$$

for some constant $c_\phi \in \mathbb{R} \setminus \{0\}$. We have $\int_0^\pi \sigma^0(\theta) \phi(\theta)^2 d\theta > 0$, according to [7, Lemm. A.2]. Hence, adjusting the constant $c_\phi$ if necessary, we can normalize $\phi$ so that $\mu \int_0^\pi \sigma^0(\theta) \phi(\theta)^2 d\theta = 1$.

**Remark 3.1.** • *For $\kappa_\sigma = \sigma_-/\sigma_+ = -1/3$, one has $0 \in \Lambda$. The singularities associated with the singular exponent 0 write $(r, \theta) \mapsto c\varphi(\theta)$ and $(r, \theta) \mapsto c \ln r \, \varphi(\theta)$, where $c$ is a constant, $\varphi(\theta) = \theta$ on $[0, \pi/4]$ and $\varphi(\theta) = (\pi - \theta)/3$ on $[\pi/4, \pi]$. As previously announced, we do not study this limit case here.*
• *For $\kappa_\sigma \in (-1, -1/3)$ such that $\kappa_\sigma \to -1^+$, there holds $\mu \to +\infty$.*
• *Finally, for $\kappa_\sigma \in \mathbb{R}^*_-$ such that $\kappa_\sigma \notin [-1, -1/3]$, there holds $\Lambda \cap \{\lambda \in \mathbb{C} \,|\, \Re e\{\lambda\} = 0\} = \emptyset$. Consequently, we can prove that $A^0$ is a Fredholm operator of index 0.*

Let $\mathscr{C}_0^\infty(\overline{\Omega} \setminus \{O\})$ refer to the set of infinitely differentiable functions supported in $\overline{\Omega} \setminus \{O\}$. Now let us give a functional framework well adapted to (4). For $\beta \in \mathbb{R}$ and $k \geq 0$, we define the Kondratiev space $V_\beta^k(\Omega)$ as the completion of $\mathscr{C}_0^\infty(\overline{\Omega} \setminus \{O\})$ for the norm

$$\|v\|_{V_\beta^k(\Omega)} := \Big( \sum_{|\alpha| \leq k} \int_\Omega r^{2(\beta+|\alpha|-k)} |\partial_{\mathbf{x}}^\alpha v|^2 \, d\mathbf{x} \Big)^{1/2}. \tag{7}$$

To take into account the Dirichlet boundary condition on $\partial\Omega$, for $\beta \in \mathbb{R}$, we introduce the space

$$\mathring{V}_\beta^1(\Omega) := \left\{ v \in V_\beta^1(\Omega) \mid v = 0 \text{ on } \partial\Omega \setminus \{O\} \right\}.$$

For the definition of the trace of the elements of $V_\beta^1(\Omega)$, we refer the reader to [27, §6.2.1]. One can check that for all $\beta \in \mathbb{R}$, $\mathring{V}_\beta^1(\Omega)$ is equal to the completion of $\mathscr{C}_0^\infty(\Omega)$ for the norm $\| \ \|_{V_\beta^1(\Omega)}$. Moreover, using a Poincaré inequality on the arc $(0,\pi)$, we can prove the estimate $\|r^{-1}v\|_\Omega \leq c \|\nabla v\|_\Omega$ for all $v \in \mathring{V}_0^1(\Omega)$. This allows to conclude that $H_0^1(\Omega) = \mathring{V}_0^1(\Omega)$. The norm in the dual space to $\mathring{V}_\beta^1(\Omega)$ is the intrinsic norm

$$\|g\|_{\mathring{V}_\beta^1(\Omega)^*} := \sup_{v \in \mathring{V}_\beta^1(\Omega) \setminus \{0\}} \frac{|\langle g, v \rangle_\Omega|}{\|v\|_{V_\beta^1(\Omega)}},$$

where $\langle \ , \ \rangle_\Omega$ refers to the duality pairing between $\mathring{V}_\beta^1(\Omega)^*$ and $\mathring{V}_\beta^1(\Omega)$. Although we adopt the same notation for the pairing between $H^{-1}(\Omega)$ and $H_0^1(\Omega)$, this will not bring further confusion. A functional framework for (4) also requires to include *a priori* information about a proper propagative behaviour of the solution at 0. For $\beta \in (0,2)$, we set

$$V_\beta^{\text{out}}(\Omega) := \text{span}\{\psi(r) r^{i\mu} \phi(\theta)\} \oplus \mathring{V}_{-\beta}^1(\Omega)$$

$$\|v\|_{V_\beta^{\text{out}}(\Omega)} := |c| + \|\tilde{v}\|_{V_{-\beta}^1(\Omega)} \quad \text{for } v(r,\theta) = c\,\psi(r) r^{i\mu} \phi(\theta) + \tilde{v}(r,\theta).$$

Here, the cut-off function $\psi \in \mathscr{C}^\infty(\mathbb{R}_+; [0,1])$ is such that $\psi(r) = 1$ for $r \leq 1$ and $\psi(r) = 0$ for $r \geq 2$. This kind of weighted spaces with detached asymptotic were introduced in [35, Chap. 6] through the generalized Green formula (see also [35, Chap. 11,12]) and the following proposition is the concretization of general results. For the proof in the particular context of this paper (transmission problems were not considered in [35]), we refer the reader to [7, Th. 4.4].

**Proposition 3.1.**
*For $\beta \in (0,2)$, define by density the operator $A_\beta^{\text{out}} : V_\beta^{\text{out}}(\Omega) \to \mathring{V}_\beta^1(\Omega)^*$ as the unique operator satisfying $\langle A_\beta^{\text{out}} u, v \rangle_\Omega = (\sigma^0 \nabla u, \nabla v)_\Omega$, $\forall u \in V_\beta^{\text{out}}(\Omega)$, $v \in \mathscr{C}_0^\infty(\Omega)$. Then the operator $A_\beta^{\text{out}}$ is Fredholm of index 0.*

It can occur, depending on the geometry and the value of the coefficient $\sigma$, that $A_\beta^{\text{out}}$ gets a non trivial kernel (of course, of finite dimension since this operator is Fredholm). We discard this possibility, considering an additional assumption

**Assumption 1.** *For $\beta \in (0,2)$, the operator $A_\beta^{\text{out}}$ is one-to-one.*

As a consequence, $A_\beta^{\text{out}}$ is an isomorphism and there exists $C > 0$ (that depends on $\beta$) such that

$$|c| + \|\tilde{v}\|_{V_{-\beta}^1(\Omega)} \leq C \, \|A_\beta^{\text{out}} v\|_{\mathring{V}_\beta^1(\Omega)^*}, \qquad \forall v = c\,\psi(r) r^{i\mu} \phi(\theta) + \tilde{v} \in V_\beta^{\text{out}}(\Omega). \tag{8}$$

**Remark 3.2.** *Assumption 1 can be slightly weaken. Indeed, in [7], it is proven that the functions of $\ker A_\beta^{\text{out}}$ belong to $\mathring{V}_{-\beta}^1(\Omega)$. In other words, they do not decompose on the strongly oscillating singularity. Thus, Assumption 1 is equivalent to assume that the only element $u \in \mathring{V}_{-\beta}^1(\Omega)$ (the important point here is that $\mathring{V}_{-\beta}^1(\Omega) \subset V_\beta^{\text{out}}(\Omega)$) such that $(\sigma \nabla u, \nabla v)_\Omega = 0$, $\forall v \in \mathring{V}_\beta^1(\Omega)$, is the null function.*

**Remark 3.3.** *Using [7, Lemm. 4.1], we can establish that Assumption 1 holds for the geometry (see Fig. 6, on right) considered in the numerical experiments of Section 7.*

## 3.2 Inner limit geometry

The second limit geometry we want to consider is obtained by applying the homothety $\mathbf{x} \mapsto \mathbf{x}/\delta$ centered at $O$. When this transformation is applied to $\Omega$, we obtain a rescaled geometry $\Xi^\delta$ where the rounded corner is fixed. Passing formally to the limit as $\delta \to 0$, the domain $\Xi^\delta$ becomes the limit domain $\Xi = \mathbb{R} \times (0, +\infty)$, see Fig. 2. The polar coordinates in the *inner limit geometry* $\Xi$ will be denoted $(\rho, \theta)$. Besides, we set $\sigma^\infty(\mathbf{x}) = \sigma_\pm$ for $\mathbf{x} \in \Xi_\pm$.

To deal with a near field problem set in this second inner limit geometry, we need to consider another, more appropriate, functional setting. For $\beta \in \mathbb{R}$, $k \geq 0$, we introduce the space $\mathcal{V}_\beta^k(\Xi)$ defined as the completion of the set $\{v|_\Xi,\, v \in \mathscr{C}_0^\infty(\mathbb{R}^2)\}$ for the norm

$$\|v\|_{\mathcal{V}_\beta^k(\Xi)} := \Big( \sum_{|\alpha| \leq k} \int_\Xi (1+\rho)^{2(\beta + |\alpha| - k)} |\partial_\mathbf{x}^\alpha v|^2\, d\mathbf{x} \Big)^{1/2}. \tag{9}$$

Introducing these spaces, our goal is to discriminate behaviours of functions at infinity only and not at $O$. This explains why the weight $\rho^{2(\beta+|\alpha|-k)}$ in (7) has been replaced by $(1+\rho)^{2(\beta+|\alpha|-k)}$ in (9) (notice also that the notation has been modified). In the sequel, we will impose the Dirichlet boundary condition on $\partial\Xi$ using the definition of the trace as in [27, §6.1.1] and working with functions belonging to

$$\mathring{\mathcal{V}}_\beta^1(\Xi) := \left\{ v \in \mathcal{V}_\beta^1(\Xi) \mid v = 0 \text{ on } \partial\Xi \right\}.$$

For all $\beta \in \mathbb{R}$, this space coincides with the completion of $\mathscr{C}_0^\infty(\Xi)$ for the norm $\|\ \|_{\mathcal{V}_\beta^1(\Xi)}$. The norm in the dual space to $\mathring{\mathcal{V}}_\beta^1(\Xi)$ is defined in a usual manner, as follows

$$\|g\|_{\mathring{\mathcal{V}}_\beta^1(\Xi)^*} := \sup_{v \in \mathring{\mathcal{V}}_\beta^1(\Xi) \setminus \{0\}} \frac{|\langle g, v\rangle_\Xi|}{\|v\|_{\mathcal{V}_\beta^1(\Xi)}},$$

where $\langle\ ,\ \rangle_\Xi$ refers to the duality pairing between $\mathring{\mathcal{V}}_\beta^1(\Xi)^*$ and $\mathring{\mathcal{V}}_\beta^1(\Xi)$. We also consider a space prescribing a proper propagative behaviour at infinity. For $\beta \in (0,2)$ we consider

$$\mathcal{V}_{-\beta}^{\text{in}}(\Xi) := \text{span}\{\chi(\rho)\rho^{-i\mu}\phi(\theta)\} \oplus \mathring{\mathcal{V}}_\beta^1(\Xi)$$

$$\|v\|_{\mathcal{V}_{-\beta}^{\text{in}}(\Xi)} := |c| + \|\tilde{v}\|_{\mathcal{V}_\beta^1(\Xi)} \quad \text{for } v(\rho,\theta) = c\,\chi(\rho)\rho^{-i\mu}\phi(\theta) + \tilde{v}(\rho,\theta).$$

Here, we take $\chi = 1 - \psi \in \mathscr{C}^\infty(\mathbb{R}_+; [0,1])$, $\psi$ being the cut-off function defined in the previous paragraph. In particular, there holds $\chi(r) = 1$ for $r \geq 2$ and $\chi(r) = 0$ for $r \leq 1$. In this second limit geometry, we have a result similar to Proposition 3.1. We do not give the proof which is very close to the one of Proposition 3.1.

**Proposition 3.2.**
*For $\beta \in (0,2)$, define by density the operator $\mathcal{A}_{-\beta}^{\text{in}} : \mathcal{V}_{-\beta}^{\text{in}}(\Xi) \to \mathring{\mathcal{V}}_{-\beta}^1(\Xi)^*$ as the unique operator satisfying $\langle \mathcal{A}_{-\beta}^{\text{in}} u, v\rangle_\Xi = (\sigma^\infty \nabla u, \nabla v)_\Xi$, $\forall u \in \mathcal{V}_{-\beta}^{\text{in}}(\Xi)$, $v \in \mathscr{C}_0^\infty(\Xi)$. Then the operator $\mathcal{A}_{-\beta}^{\text{in}}$ is Fredholm of index 0.*

In the present case, we need to impose injectivity of $\mathcal{A}_{-\beta}^{\text{in}}$.

**Assumption 2.** *For $\beta \in (0,2)$, the operator $\mathcal{A}_{-\beta}^{\text{in}}$ is one-to-one.*

Thus, $\mathcal{A}_{-\beta}^{\text{in}}$ is an isomorphism and there exists $C > 0$ (that depends on $\beta$) such that

$$|c| + \|\tilde{v}\|_{\mathcal{V}_\beta^1(\Xi)} \leq C\, \|\mathcal{A}_{-\beta}^{\text{in}} v\|_{\mathring{\mathcal{V}}_{-\beta}^1(\Xi)^*}, \qquad \forall v = c\,\chi(\rho)\rho^{-i\mu}\phi(\theta) + \tilde{v} \in \mathcal{V}_{-\beta}^{\text{in}}(\Xi)\,. \tag{10}$$

**Remark 3.4.** *As for Remark 3.2, we can prove that the functions of $\text{Ker}\,\mathcal{A}_{-\beta}^{\text{in}}$ belong to $\mathring{\mathcal{V}}_\beta^1(\Xi)$. Therefore, Assumption 2 is equivalent to assume that the only element $u \in \mathring{\mathcal{V}}_\beta^1(\Xi)$ such that $(\sigma^\infty \nabla u, \nabla v)_\Xi = 0$, $\forall v \in \mathring{\mathcal{V}}_{-\beta}^1(\Xi)$, is the null function. Moreover, [7, Lemm. 4.1] allows to show that Assumption 2 holds for the configuration investigated in Section 7.*

# 4 Asymptotic description of the solution

In this section, we apply the matched expansions technique [51, 24] in order to derive the first orders of the expansion of $u^\delta$. In the present problem, there appears a boundary layer in the neighbourhood of the rounded corner, so we have to compute an inner and an outer expansion. As tools for separating outer and inner regions, we reintroduce the $\mathscr{C}^\infty$ cut-off functions $\chi$ and $\psi$ such that
$$\chi(r) + \psi(r) = 1, \qquad \chi(r) = 1 \text{ for } r \geq 2, \quad \text{and} \quad \chi(r) = 0 \text{ for } r \leq 1.$$

In the sequel, we shall denote $\chi_t(r) = \chi(r/t)$ and $\psi_t(r) = \psi(r/t)$ for $t > 0$. We shall also set $(\tau_t \cdot v)(r, \theta) := v(r/t, \theta)$ for $v \in L^2(\Omega)$. Then we decompose the right hand side $f \in H^{-1}(\Omega)$ into an outer and an inner contribution,

$$f = f_{\text{ext}}^\delta + \tau_\delta \cdot F_{\text{in}}^\delta \qquad \text{with} \qquad \left| \begin{array}{l} \langle f_{\text{ext}}^\delta, v \rangle_\Omega := \langle f, v\, \chi_{\sqrt{\delta}} \rangle_\Omega \\ \langle F_{\text{in}}^\delta, v \rangle_\Xi := \langle f, (\tau_\delta \cdot v)\, \psi_{\sqrt{\delta}} \rangle_\Omega \\ \langle \tau_\delta \cdot F_{\text{in}}^\delta, v \rangle_\Omega := \langle F_{\text{in}}^\delta, \tau_{1/\delta} \cdot v \rangle_\Xi \end{array} \right. \tag{11}$$

In the case where $f \in L^2(\Omega)$, this decomposition simply reads $f(r, \theta) = f_{\text{ext}}^\delta(r, \theta) + \delta^{-2} F_{\text{in}}^\delta(r/\delta, \theta)$ with $f_{\text{ext}}^\delta(r, \theta) = f(r, \theta)\chi(r/\sqrt{\delta})$, and $F_{\text{in}}^\delta(\rho, \theta) = \delta^2 f(\delta\rho, \theta)\psi(\rho\sqrt{\delta})$. Observe that $\text{supp}(f_{\text{ext}}^\delta)$ excludes $D(O, \sqrt{\delta})$ while $\text{supp}(F_{\text{in}}^\delta)$ is included in the closure of $D(O, 2/\sqrt{\delta})$. As a consequence, we have
$$f_{\text{ext}}^\delta \in \mathring{V}_\beta^1(\Omega)^* \quad \text{and} \quad F_{\text{in}}^\delta \in \mathring{\mathcal{V}}_{-\beta}^1(\Xi)^*, \quad \forall \beta \in \mathbb{R}.$$

## 4.1 Outer expansion

In the outer region, we look for an expansion of $u^\delta(r, \theta)$ in the form
$$u^\delta(r, \theta) = u_0^\delta(r, \theta) + a(\delta)\, \zeta(r, \theta) + \cdots, \tag{12}$$

where $u_0^\delta$, $\zeta$ and $a(\delta)$ have to be defined. As regards the first term of this expansion, plugging (12) into (1), we see that the following equation has to hold
$$\left| \begin{array}{l} u_0^\delta \in V_\beta^{\text{out}}(\Omega) \text{ such that} \\ -\text{div}(\sigma^0 \nabla u_0^\delta) = f_{\text{ext}}^\delta \quad \text{in } \Omega. \end{array} \right. \tag{13}$$

According to Assumption 1, the equation above uniquely defines $u_0^\delta$. Moreover, by definition of $V_\beta^{\text{out}}(\Omega)$, we have
$$u_0^\delta(r, \theta) = c_0^\delta\, \psi(r) r^{i\mu}\, \phi(\theta) + \tilde{u}_0^\delta(r, \theta) \qquad \text{where} \quad c_0^\delta \in \mathbb{C} \quad \text{and} \quad \tilde{u}_0^\delta \in \mathring{V}_{-\beta}^1(\Omega). \tag{14}$$

Concerning the second term on the right hand side of (12), we still do not know what is a relevant form for the jauge function $a(\delta)$, we shall determine this when applying the matching procedure later on. Plugging (12) into (1), we find
$$\text{div}(\sigma^0 \nabla \zeta) = 0 \quad \text{in } \Omega, \quad \zeta = 0 \text{ on } \partial\Omega. \tag{15}$$

We have to look for $\zeta$ in a space wider than just $V_\beta^{\text{out}}(\Omega)$, otherwise we would have had to conclude that $\zeta = 0$ which is excluded. For $\beta \in \mathbb{R}$, let us consider $A_\beta : \mathring{V}_\beta^1(\Omega) \to \mathring{V}_{-\beta}^1(\Omega)^*$ defined as the unique operator satisfying $\langle A_\beta u, v \rangle_\Omega = (\sigma\nabla u, \nabla v)_\Omega$, $\forall u \in \mathring{V}_\beta^1(\Omega), v \in \mathring{V}_{-\beta}^1(\Omega)$. For $\beta = 0$, since $\mathring{V}_0^1(\Omega) = H_0^1(\Omega)$, there holds $A_0 = A^0$, where $A^0$ has been defined in §2.2. Equation (15) imposes that $\zeta \in \text{Ker}\, A_\beta$. Choosing $\beta \in (0, 2)$, we know $\zeta$ completely, up to some multiplicative constant, according to the following proposition.

**Proposition 4.1.**
*For $\beta \in (0, 2)$, under Assumption 1, we have* $\dim(\text{Ker}\, A_\beta) = 1$.

**Proof:** This result is a consequence of the general theorem on index [35, Chap. 4, Th. 3.3]. For ease of reading, we provide the proof for this particular problem. Set $w_1(r,\theta) = \psi(r)r^{-i\mu}\phi(\theta)$. According to the construction of $\phi(\theta)$, we have $g = \text{div}(\sigma^0 \nabla w_1) \in \mathscr{C}_0^\infty(\Omega)$. Hence, according to Proposition 3.1 and Assumption 1, there exists a unique $w_2 \in V_\beta^{\text{out}}(\Omega)$ such that $\text{div}(\sigma^0 \nabla w_2) = g$ in $\Omega$, so that $w := w_1 - w_2 \in \text{Ker}\, A_\beta$. Notice that $w \neq 0$ since $w - w_1 \in V_\beta^{\text{out}}(\Omega)$ and $w_1 \notin V_\beta^{\text{out}}(\Omega)$. Thus, $\dim(\text{Ker}\, A_\beta) \geq 1$.

Assume that $v \in \mathring{V}_\beta^1(\Omega)$ is another element of $\text{Ker}\, A_\beta$. Since $\beta \in (0,2)$ and $\Lambda \cap \{\lambda \in \mathbb{C} \mid |\Re e\{\lambda\}| \leq \beta\} = \{+i\mu, -i\mu\}$, according to [7, Th. 5.2], there exists $c_-, c_+ \in \mathbb{C}$ such that $v(r,\theta) - (c_+\psi(r)r^{i\mu} + c_-\psi(r)r^{-i\mu})\phi(\theta) \in \mathring{V}_{-\beta}^1(\Omega)$. As a consequence, $v - c_- w \in V_\beta^{\text{out}}(\Omega)$, and $v - c_- w \in \text{Ker}\, A_\beta^{\text{out}}$, which implies $v - c_- w = 0$ according to Proposition 3.1 and Assumption 1. Thus, $\dim(\text{Ker}\, A_\beta) \leq 1$. $\square$

Equation (15) imposes that $\zeta \in \text{Ker}\, A_\beta$. Adjusting the jauge function $a(\delta)$ if necessary, we can take $\zeta$ as the only element of $\text{Ker}\, A_\beta$ admitting the following expansion

$$\zeta(r,\theta) = \psi(r)r^{-i\mu}\phi(\theta) + c_\zeta\,\psi(r)r^{i\mu}\phi(\theta) + \tilde{\zeta}(r,\theta)$$
$$\text{with } \tilde{\zeta} \in \mathring{V}_{-\beta}^1(\Omega) \text{ and } \beta \in (0,2)\,. \tag{16}$$

**Lemma 4.1.**
*The constant $c_\zeta \in \mathbb{C}$ coming into play in the radial expansion (16) verifies $|c_\zeta| = 1$.*

**Proof:** To prove this result, let us apply the method of [32] (see also [35, Chap. 6]). Since $\zeta \in \text{Ker}\, A_\beta$, according to the definition of $A_\beta$, we have $0 = \Im m\{\langle A_\beta \zeta, \chi_\varepsilon \zeta\rangle_\Omega\} = \Im m\{\int_\Omega \sigma^0 \nabla \zeta \cdot \nabla(\chi_\varepsilon \overline{\zeta})d\mathbf{x}\} = \Im m\{\int_\Omega \sigma^0 \overline{\zeta} \partial_r \zeta\, \partial_r \chi_\varepsilon d\mathbf{x}\}$. We take a look at the behaviour of this expression when $\varepsilon \to 0$. Using (16), we have

$$\begin{aligned}
0 = \Im m\{\int_\Omega \sigma^0 \overline{\zeta}\partial_r\zeta\,\partial_r\chi_\varepsilon d\mathbf{x}\} &= \int_0^\pi \int_\varepsilon^{2\varepsilon} \sigma^0(\theta)\,\Im m\{\overline{\zeta}\partial_r\zeta\}\partial_r\chi_\varepsilon\,rdrd\theta \\
&= \mu \int_0^\pi \int_\varepsilon^{2\varepsilon} \sigma^0(\theta)\,|\phi(\theta)|^2(|c_\zeta|^2 - 1)\partial_r\chi_\varepsilon\,drd\theta + O(\varepsilon^\beta) \\
&= (|c_\zeta|^2 - 1)\,\mu \int_0^\pi \sigma^0(\theta)\,|\phi(\theta)|^2 d\theta \int_\varepsilon^{2\varepsilon} \partial_r\chi_\varepsilon\,dr + O(\varepsilon^\beta) \\
&= (|c_\zeta|^2 - 1) + O(\varepsilon^\beta).
\end{aligned}$$

This implies $|c_\zeta|^2 - 1 = 0$. $\square$

To sum up, from a formal point of view, up to some remainder with respect to $\delta$, we have the following asymptotic behaviour for $r \to 0$,

$$\begin{aligned}
u^\delta(r,\theta) &= u_0^\delta(r,\theta) + a(\delta)\,\zeta(r,\theta) + \cdots \\
&= (c_0^\delta + a(\delta)c_\zeta)\,r^{i\mu}\phi(\theta) + a(\delta)\,r^{-i\mu}\phi(\theta) + \cdots.
\end{aligned} \tag{17}$$

### 4.2 Inner expansion

In the inner region, close to $O$, we consider the change of coordinates $r = \delta\rho$, we set $U^\delta(\rho,\theta) = u^\delta(r,\theta)$, and look for an expansion of $U^\delta(\rho,\theta)$ of the form

$$U^\delta(\rho,\theta) = U_0^\delta(\rho,\theta) + A(\delta)\,Z(\rho,\theta) + \cdots, \tag{18}$$

where $U_0^\delta$, $Z$ and $A(\delta)$ have to be defined. Starting from (1), making the change of variables $\rho = \delta^{-1}r$, we formally obtain

$$\left| \begin{array}{l} U_0^\delta \in \mathcal{V}_{-\beta}^{\text{in}}(\Xi) \text{ such that} \\ -\text{div}(\sigma^\infty \nabla U_0^\delta) = F_{\text{in}}^\delta. \end{array} \right. \tag{19}$$

According to Proposition 3.2, $U_0^\delta$ is uniquely determined. Note that, by the very definition of $\mathcal{V}_{-\beta}^{\mathrm{in}}(\Xi)$, we have the following expansion

$$U_0^\delta(\rho,\theta) = C_0^\delta\,\chi(\rho)\rho^{-i\mu}\,\phi(\theta) + \tilde{U}_0^\delta(\rho,\theta) \qquad \text{where} \quad C_0^\delta \in \mathbb{C} \quad \text{and} \quad \tilde{U}_0^\delta \in \mathring{\mathcal{V}}_\beta^1(\Xi). \tag{20}$$

Now, let us take a look at the second term in the asymptotics of $U^\delta$. Like for the outer expansion, we still have not determined a relevant form for the jauge function $A(\delta)$. This will be a result of the matching procedure. Plugging (18) into (1) yields

$$\operatorname{div}(\sigma^\infty \nabla Z) = 0 \quad \text{in } \Xi, \quad Z = 0 \ \text{on } \partial\Xi\,. \tag{21}$$

The function $Z$ has to be looked for in a space wider than $\mathcal{V}_{-\beta}^{\mathrm{in}}(\Xi)$, otherwise, we would have to conclude $Z = 0$, which is not relevant. For $\beta \in \mathbb{R}$, define $\mathcal{A}_{-\beta} : \mathring{\mathcal{V}}_{-\beta}^1(\Xi) \to \mathring{\mathcal{V}}_\beta^1(\Xi)^*$ the unique operator satisfying $\langle \mathcal{A}_{-\beta}u,v\rangle_\Xi = (\sigma^\infty \nabla u, \nabla v)_\Xi,\ \forall u \in \mathring{\mathcal{V}}_{-\beta}^1(\Xi),\ v \in \mathring{\mathcal{V}}_\beta^1(\Xi)$. Equations (21) determine $Z \in \mathring{\mathcal{V}}_{-\beta}^1(\Xi)$, for $\beta \in (0,2)$, up a multiplicative constant, according to the following proposition.

**Proposition 4.2.**
*For $\beta \in (0,2)$, under Assumption 2, we have $\dim(\operatorname{Ker}\mathcal{A}_{-\beta}) = 1$.*

The proof of this Proposition follows the same lines as for Proposition 4.1, hence we do not present it. Adjusting the jauge function $A(\delta)$ if necessary, we can take $Z(\rho,\theta)$ as the unique element of $\operatorname{Ker}\mathcal{A}_{-\beta}$ that admits the expansion

$$\begin{gathered}
Z(\rho,\theta) = \chi(\rho)\rho^{i\mu}\phi(\theta) + C_{\mathrm{z}}\,\chi(\rho)\rho^{-i\mu}\,\phi(\theta) + \tilde{Z}(\rho,\theta) \\
\text{with } \tilde{Z} \in \mathring{\mathcal{V}}_\beta^1(\Xi) \text{ and } \beta \in (0,2)\,.
\end{gathered} \tag{22}$$

The existence of such an expansion for $Z(\rho,\theta)$ is a consequence of a classical residue calculus involving the Mellin symbol (with respect to $\rho$) of $\mathcal{A}_{-\beta}$. Again, we use the fact that, for $\beta \in (0,2)$, $\Lambda \cap \{\lambda \in \mathbb{C} \mid |\Re e\{\lambda\}| \leq \beta\} = \{+i\mu, -i\mu\}$, where $\Lambda$ is the set of singular exponents defined in (5). Concerning the coefficient $C_{\mathrm{z}}$ we have a result similar to Lemma 4.1.

**Lemma 4.2.**
*The constant $C_{\mathrm{z}} \in \mathbb{C}$ coming into play in the radial expansion (22) verifies $|C_{\mathrm{z}}| = 1$.*

To conclude, up to some remainder with respect to $\delta$, the field in the inner region admits the following behaviour for $\rho = r/\delta \to +\infty$,

$$\begin{aligned}
U^\delta(\rho,\theta) &= U_0^\delta(\rho,\theta) + A(\delta)\,Z(\rho,\theta) + \cdots \\
&= A(\delta)\rho^{i\mu}\phi(\theta) + (C_0^\delta + A(\delta)C_{\mathrm{z}})\rho^{-i\mu}\phi(\theta) + \cdots,
\end{aligned}$$

$$\text{i.e.,} \quad u^\delta(r,\theta) = A(\delta)\Big(\frac{r}{\delta}\Big)^{i\mu}\phi(\theta) + (C_0^\delta + A(\delta)C_{\mathrm{z}})\Big(\frac{r}{\delta}\Big)^{-i\mu}\phi(\theta) + \cdots. \tag{23}$$

### 4.3 The matching procedure

In order to conclude the construction of the first terms of the asymptotics of $u_\delta$, we still have to determine the jauge functions $a(\delta)$ and $A(\delta)$. This will be achieved by applying the matching procedure to the far field and near field expansions. Matching procedure has been described in detail in the reference books [51, 24] and [31, Chap. 2]. In the present case it consists in equating expansions (17) and (23) for $r, \delta \to 0$,

$$(c_0^\delta + a(\delta)c_\zeta)\,r^{i\mu}\phi(\theta) + a(\delta)\,r^{-i\mu}\phi(\theta) = A(\delta)\Big(\frac{r}{\delta}\Big)^{i\mu}\phi(\theta) + (C_0^\delta + A(\delta)C_{\mathrm{z}})\Big(\frac{r}{\delta}\Big)^{-i\mu}\phi(\theta). \tag{24}$$

As $r^{i\mu}$ and $r^{-i\mu}$ are independent, Equation (24) leads to the following matching conditions

$$c_0^\delta + a(\delta)c_\zeta = A(\delta)\,\delta^{-i\mu} \qquad \text{and} \qquad a(\delta) = (C_0^\delta + A(\delta)C_z)\,\delta^{i\mu}. \tag{25}$$

In these equations, $c_0^\delta, c_\zeta, C_0^\delta, C_z$ are known data. Equations (25) are uniquely solvable with respect to $a(\delta), A(\delta)$ under the following condition

$$\delta^{-2i\mu} \neq c_\zeta\, C_z. \tag{26}$$

This condition requires comments. First, notice that the matching of expansions is subordinated to a condition involving $\delta$. Observe also that, since $|c_\zeta C_z| = 1$ according to Lemmas 4.1 and 4.2, there exists $\delta_\star \in (0, +\infty)$ such that $c_\zeta\, C_z = \delta_\star^{-2i\mu}$. Hence the condition (26) is violated whenever there exists some $k \in \mathbb{Z}$ such that $\delta = (e^{\pi/\mu})^k \delta_\star$. This is an uncomfortable situation because the set $I_\star := \{\,(e^{\pi/\mu})^k \delta_\star, k \in \mathbb{Z}\,\}$ admits 0 as accumulation point, so that (26) cannot be garantied simply by imposing a condition of the form "$\delta \in (0, \delta_0]$ for some $\delta_0 > 0$ small enough". For further investigations concerning this type of phenomena in matched asymptotic methods, we refer the reader to [24, 30, 15] and to [31, Chap. 4]. In the case where (26) holds, we have the following expression for the jauge functions

$$a(\delta) = \frac{c_0^\delta\, C_z\, \delta^{2i\mu} + C_0^\delta\, \delta^{i\mu}}{1 - (\delta/\delta_\star)^{2i\mu}} \qquad \text{and} \qquad A(\delta) = \frac{C_0^\delta\, c_\zeta\, \delta^{2i\mu} + c_0^\delta\, \delta^{i\mu}}{1 - (\delta/\delta_\star)^{2i\mu}}. \tag{27}$$

Because of these relations, as $\delta \to 0$, both jauge functions are *a priori* unbounded, which discards any perspective of establishing error estimates. For this reason, we have to introduce restrictions on the set where $\delta$ varies. In particular, we have to exclude a neighbourhood of $I_\star$.

**Definition 4.1.** *For $\alpha \in (0, 1/2)$, define* $I(\alpha) := \bigcup_{k=-\infty}^{+\infty} \left[\,\delta_\star\, e^{(k+\alpha)\frac{\pi}{\mu}},\ \delta_\star\, e^{(k+1-\alpha)\frac{\pi}{\mu}}\,\right]$.

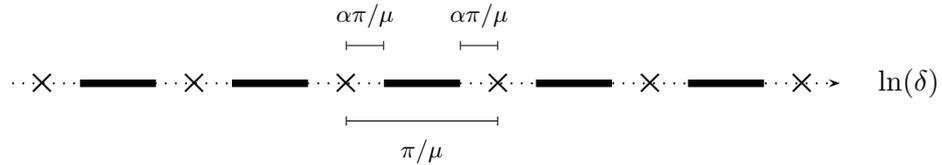

Fig. 5: The crosses represent the set $\{\ln(\delta), \delta \in I_\star\}$. The thick segments represent the set $\{\ln(\delta), \delta \in I(\alpha)\}$.

An elementary calculus shows that $I(\alpha) = \{\delta \in (0, +\infty)\,|\,1/(1 - (\delta/\delta_\star)^{2i\mu}) \leq 1/(2\sin\alpha)\}$. Therefore, as illustrated in Fig. 5, there holds $I(\alpha) \cap I_\star = \emptyset$. Note also that $I(\alpha)$ admits 0 as accumulation point. In the remaining of this paper, the statement "$\delta \to 0$" should be understood in the sense that $\delta \to 0$ *and* $\delta \in I(\alpha)$ for some $\alpha \in (0, 1/2)$. This assumption implies that the jauge functions remain bounded as $\delta \to 0$, *i.e.* there exists $c_\alpha > 0$ independent of $\delta$ such that

$$|a(\delta)| + |A(\delta)| \leq c_\alpha\,(\,|c_0^\delta| + |C_0^\delta|\,), \qquad \forall \delta \in I(\alpha). \tag{28}$$

## 5 Construction and norm estimate of the resolvent

The main purpose of the present section is Theorem 5.1 that provides an estimate of $(A^\delta)^{-1}$ with respect to $\delta \in I(\alpha)$ in a well chosen norm. In a first step, using the asymptotic description of the previous section, we will build explicitly the inverse of $A^\delta$. Then, the proof of Theorem 5.1 will be the result of two key estimates provided by Proposition 5.1 and 5.2. Note that, in the sequel, $c > 0$ shall denote a constant independent of $\delta$ that may vary from one line to another.

## 5.1 Explicit expression of the resolvent

Thanks to a trick of [31, Chap. 2], [34] which relies on the use of overlapping cut-off functions, we construct an approximation $\hat{u}^\delta$ of $u_\delta$ by means of the outer and the inner expansions. First of all, going back to the definition of the cut-off functions $\chi_\delta$, $\psi$ (see the beginning of §4), observe that $\chi_\delta(r) + \psi(r) - \chi_\delta(r)\psi(r) = 1$. As an approximation of the field $u_\delta$ from the inner and outer expansion, we set

$$\hat{u}^\delta(r,\theta) = u^\delta_{\text{ext}}(r,\theta)\chi_\delta(r) + U^\delta_{\text{in}}(r/\delta,\theta)\psi(r) - m^\delta(r,\theta)\chi_\delta(r)\psi(r)$$

$$\text{with} \quad \begin{vmatrix} u^\delta_{\text{ext}}(r,\theta) &=& u^\delta_0(r,\theta) + a(\delta)\zeta(r,\theta), \\ U^\delta_{\text{in}}(\rho,\theta) &=& U^\delta_0(\rho,\theta) + A(\delta)Z(\rho,\theta), \\ m^\delta(r,\theta) &=& A(\delta)(r/\delta)^{i\mu}\phi(\theta) + a(\delta)r^{-i\mu}\phi(\theta)\,. \end{vmatrix} \quad (29)$$

We recall that we have introduced the map $\tau_\delta$, such that, for any function $v$, $(\tau_\delta \cdot v)(r,\theta) := v(r/\delta,\theta)$. Let us check that $\hat{u}^\delta$ satisfies Equation (1) up to some remainder. First of all, it is clear, from the definition of $u^\delta_0, U^\delta_0, \zeta, Z$ and $\phi(\theta)$, that $\hat{u}^\delta = 0$ on $\partial\Omega$ for $\delta$ small enough. For any $v \in \mathrm{H}^1_0(\Omega)$, according to the definition of the operators $\mathrm{A}^{\text{out}}_\beta$ and $\mathcal{A}^{\text{in}}_{-\beta}$ (see Propositions 3.1 and 3.2), we have

$$\begin{aligned}(\sigma^\delta \nabla \hat{u}^\delta, \nabla v)_\Omega &=& (\sigma^0 \nabla(u^\delta_{\text{ext}}\chi_\delta), \nabla v)_\Omega + (\sigma^\delta \nabla(\tau_\delta \cdot (U^\delta_{\text{in}}\psi_{1/\delta})), \nabla v)_\Omega - (\sigma^\delta \nabla(m^\delta \chi_\delta \psi), \nabla v)_\Omega \\ &=& \langle f^\delta_{\text{ext}}, \chi_\delta v \rangle_\Omega + (\sigma^0 \nabla\chi_\delta, \overline{u}^\delta_{\text{ext}}\nabla v - v\nabla \overline{u}^\delta_{\text{ext}})_\Omega \\ && + \langle F^\delta_{\text{in}}, \psi_{1/\delta}(\tau_{1/\delta}\cdot v)\rangle_\Xi + (\sigma^\infty \nabla\psi_{1/\delta}, \overline{U}^\delta_{\text{in}}\nabla(\tau_{1/\delta}\cdot v) - (\tau_{1/\delta}\cdot v)\nabla \overline{U}^\delta_{\text{in}})_\Xi \\ && - (\sigma^\delta \nabla(\chi_\delta \psi), \overline{m}^\delta \nabla v - v\nabla \overline{m}^\delta)_\Omega. \end{aligned}$$

In the calculus above, we used the fact that $\sigma^\delta = \sigma^0$ on $\text{supp}(\chi_\delta)$ and $\sigma^\delta = \tau_\delta \cdot \sigma^\infty$ on $\text{supp}(\psi)$. We also used that $(\sigma^\delta \nabla m^\delta, \nabla v)_\Omega = 0$ for any $v \in \mathrm{H}^1_0(\Omega)$ such that $\text{supp}(v) \subset \text{supp}(\chi_\delta \psi)$ according to the construction of $\phi(\theta)$. Now, observe that $f^\delta_{\text{ext}} \chi_\delta = f^\delta_{\text{ext}}$ and $F^\delta_{\text{in}} \psi_{1/\delta} = F^\delta_{\text{in}}$ for $\delta$ small enough. Recall that $f = f^\delta_{\text{ext}} + \tau_\delta \cdot F^\delta_{\text{in}}$, and $-\chi_\delta \psi = 1 - \chi_\delta - \psi$. Let us denote $M^\delta := \tau_{1/\delta}\cdot m^\delta$. We obtain

$$\begin{aligned}(\sigma^\delta \nabla \hat{u}^\delta, \nabla v)_\Omega &=& \langle f, v\rangle_\Omega \\ && + (\sigma^0 \nabla\chi_\delta, (\overline{u}^\delta_{\text{ext}} - \overline{m}^\delta)\nabla v - v\nabla(\overline{u}^\delta_{\text{ext}} - \overline{m}^\delta))_\Omega \\ && + (\sigma^\infty \nabla\psi_{1/\delta}, (\overline{U}^\delta_{\text{in}} - \overline{M}^\delta)\nabla(\tau_{1/\delta}\cdot v) - (\tau_{1/\delta}\cdot v)\nabla(\overline{U}^\delta_{\text{in}} - \overline{M}^\delta))_\Xi. \end{aligned} \quad (30)$$

The calculus above shows that $\hat{u}^\delta$ satisfies the same equation as (1), up to some remainder. Note that the dependency of $u^\delta_0, U^\delta_0$ and $a(\delta), A(\delta)$ is linear with respect to $f$. Hence, we can consider operators $\hat{\mathrm{R}}^\delta : \mathrm{H}^{-1}(\Omega) \to \mathrm{H}^1_0(\Omega)$ and $\mathrm{K}^\delta : \mathrm{H}^{-1}(\Omega) \to \mathrm{H}^{-1}(\Omega)$ respectively defined by

$$\begin{aligned}\hat{\mathrm{R}}^\delta f &= \hat{u}^\delta \\ \langle \mathrm{K}^\delta f, v\rangle_\Omega &= (\sigma^0 \nabla\chi_\delta, (\overline{u}^\delta_{\text{ext}} - \overline{m}^\delta)\nabla v - v\nabla(\overline{u}^\delta_{\text{ext}} - \overline{m}^\delta))_\Omega \\ &\quad + (\sigma^\infty \nabla\psi_{1/\delta}, (\overline{U}^\delta_{\text{in}} - \overline{M}^\delta)\nabla(\tau_{1/\delta}\cdot v) - (\tau_{1/\delta}\cdot v)\nabla(\overline{U}^\delta_{\text{in}} - \overline{M}^\delta))_\Xi \end{aligned} \quad (31)$$

for any $f \in \mathrm{H}^{-1}(\Omega)$ and $v \in \mathrm{H}^1_0(\Omega)$. According to (30), we have the simple identity

$$\mathrm{A}^\delta \cdot \hat{\mathrm{R}}^\delta = \mathrm{Id} + \mathrm{K}^\delta, \quad (32)$$

where we recall that $\mathrm{A}^\delta : \mathrm{H}^1_0(\Omega) \to \mathrm{H}^{-1}(\Omega)$ is the operator such that $\langle \mathrm{A}^\delta u, v\rangle_\Omega = (\sigma^\delta \nabla u, \nabla v)_\Omega$, $\forall u, v \in \mathrm{H}^1_0(\Omega)$. We are going to prove that $\mathrm{K}^\delta$ is "small" so that $\hat{\mathrm{R}}^\delta \cdot (\mathrm{Id} + \mathrm{K}^\delta)^{-1}$ is a right inverse of $\mathrm{A}^\delta$. Since $\mathrm{A}^\delta$ is a self-adjoint operator, this will prove that $\hat{\mathrm{R}}^\delta \cdot (\mathrm{Id} + \mathrm{K}^\delta)^{-1}$ is the resolvent of $A^\delta$. To proceed, we need an estimate of $K^\delta$ in the operator norm for $\delta \to 0$. This estimate will be formulated in a $\delta$-dependent norm. We set

$$\begin{aligned}
\|v\|_{V^1_{\beta,\delta}(\Omega)} &:= \Big( \int_\Omega (\delta+r)^{2\beta}|\nabla v|^2 rdrd\theta + \int_\Omega (\delta+r)^{2(\beta-1)}|v|^2 rdrd\theta \Big)^{1/2}, \\
\|g\|_{\mathring{V}^1_{\beta,\delta}(\Omega)^*} &:= \sup_{v \in H^1_0(\Omega)\setminus\{0\}} \frac{|\langle g,v\rangle|}{\|v\|_{V^1_{\beta,\delta}(\Omega)}}.
\end{aligned} \tag{33}$$

For a fixed $\delta > 0$, there is a clear equivalence between $\|\ \|_{V^1_{\beta,\delta}(\Omega)}$ and the usual norm $\|\ \|_{H^1_0(\Omega)}$. However, equivalence constants depend on $\delta$ and, as $\delta \to 0$, these two norms adopt a different behaviour. The norms (33) are well adapted to asymptotic estimates in the present case, because they share similarities with genuine weighted norms.

**Proposition 5.1.**
*For any continuous operator $T : H^{-1}(\Omega) \to H^{-1}(\Omega)$, consider the following $\delta$-dependent operator norm,*

$$\|T, \mathring{V}^1_{\beta,\delta}(\Omega)^* \to \mathring{V}^1_{\beta,\delta}(\Omega)^*\| := \sup_{g \in H^{-1}(\Omega)\setminus\{0\}} \frac{\|Tg\|_{\mathring{V}^1_{\beta,\delta}(\Omega)^*}}{\|g\|_{\mathring{V}^1_{\beta,\delta}(\Omega)^*}}.$$

*Then for any $\beta \in (0,2)$, and for any $\varepsilon > 0$ such that $\beta > \varepsilon$ and $\beta < 2-\varepsilon$, there exists a constant $c > 0$ independent of $\delta$ such that*

$$\|K^\delta, \mathring{V}^1_{\beta,\delta}(\Omega)^* \to \mathring{V}^1_{\beta,\delta}(\Omega)^*\| \leq c\,(\delta^{\varepsilon/2} + \delta^{\beta-\varepsilon/2}), \qquad \forall \delta \in I(\alpha) \cap (0,1). \tag{34}$$

*In particular, the operator $\text{Id} + K^\delta$ is invertible for $\delta \in I(\alpha)$ small enough.*

**Proof:** For a given $f \in H^{-1}(\Omega)$, we bound $\|K^\delta f\|_{\mathring{V}^1_{\beta,\delta}(\Omega)^*}$ by means of $\|f\|_{\mathring{V}^1_{\beta,\delta}(\Omega)^*}$. In this optic, observe that $\sup_{\mathbb{R}} |r\partial_r \chi_\delta|$ and $\sup_{\mathbb{R}} |\rho\partial_\rho \psi_{1/\delta}|$ remain uniformly bounded as $\delta \to 0$. Moreover, $K^\delta f$ only depends on the values of $u^\delta_{\text{ext}} - m^\delta$ and $U^\delta_{\text{in}} - M^\delta$ in the regions $\mathbb{Q}^\delta$ and $\mathbb{Q}^{1/\delta}$, with $\mathbb{Q}^t := \{\mathbf{x}(r,\theta) \in \mathbb{R} \times (0,+\infty) \mid t \leq r \leq 2t\}$ for $t > 0$. Consequently, according to (31), for $\beta, \varepsilon \in \mathbb{R}$ whose values will be specified later, there holds

$$\begin{aligned}
|\langle K^\delta f, v\rangle_\Omega| &\leq c\,\|u^\delta_{\text{ext}} - m^\delta\|_{V^1_{-\beta-\varepsilon}(\mathbb{Q}^\delta)} \|v\|_{V^1_{\beta+\varepsilon}(\mathbb{Q}^\delta)} \\
&\quad + c\,\|U^\delta_{\text{in}} - M^\delta\|_{\mathcal{V}^1_{\beta-\varepsilon}(\mathbb{Q}^{1/\delta})} \|\tau_{1/\delta} \cdot v\|_{\mathcal{V}^1_{-\beta+\varepsilon}(\mathbb{Q}^{1/\delta})}.
\end{aligned} \tag{35}$$

In the estimate above we used the fact that $\text{supp}(\nabla\chi_\delta) \subset \mathbb{Q}^\delta$ and $\text{supp}(\nabla\psi_{1/\delta}) \subset \mathbb{Q}^{1/\delta}$. We first bound the terms involving the test function $v \in H^1_0(\Omega)$:

$$\begin{aligned}
\|v\|_{V^1_{\beta+\varepsilon}(\mathbb{Q}^\delta)} &\leq c\,\delta^\varepsilon \|v\|_{V^1_\beta(\mathbb{Q}^\delta)} \leq c\,\delta^\varepsilon \|v\|_{V^1_{\beta,\delta}(\mathbb{Q}^\delta)} \leq c\,\delta^\varepsilon \|v\|_{V^1_{\beta,\delta}(\Omega)} \\
\|\tau_{1/\delta} \cdot v\|_{\mathcal{V}^1_{-\beta+\varepsilon}(\mathbb{Q}^{1/\delta})} &= \delta^{\beta-\varepsilon} \|v\|_{V^1_{-\beta+\varepsilon,\delta}(\mathbb{Q}^1)} \leq c\,\delta^{\beta-\varepsilon} \|v\|_{V^1_{\beta,\delta}(\mathbb{Q}^1)} \leq c\,\delta^{\beta-\varepsilon} \|v\|_{V^1_{\beta,\delta}(\Omega)}.
\end{aligned} \tag{36}$$

Plugging (36) into (35), dividing by $\|v\|_{V^1_{\beta,\delta}(\Omega)}$ and taking the supremum over all $v \in H^1_0(\Omega)$, we obtain

$$\|K^\delta f\|_{\mathring{V}^1_{\beta,\delta}(\Omega)^*} \leq c\,\delta^\varepsilon \|u^\delta_{\text{ext}} - m^\delta\|_{V^1_{-\beta-\varepsilon}(\mathbb{Q}^\delta)} + c\,\delta^{\beta-\varepsilon} \|U^\delta_{\text{in}} - M^\delta\|_{\mathcal{V}^1_{\beta-\varepsilon}(\mathbb{Q}^{1/\delta})}. \tag{37}$$

We first bound the terms associated with $u^\delta_{\text{ext}} - m^\delta$. On $\mathbb{Q}^\delta$, a direct computation using (27) yields $u^\delta_{\text{ext}} - m^\delta = \tilde{u}^\delta_0 + a(\delta)\tilde{\zeta} \in \mathring{V}^1_\gamma(\Omega)$ for all $\gamma \in (-2,0)$ (we recall that $\tilde{u}^\delta_0$ and $\tilde{\zeta}$ are respectively defined in (14) and (16)). Set $\beta \in (0,2)$. Choose $\varepsilon > 0$ such that $\beta + \varepsilon < 2$ and $\beta - \varepsilon > 0$. According to (28), (8), (13), (10) and (19), one has

$$\begin{aligned}
\|u^\delta_{\text{ext}} - m^\delta\|_{V^1_{-\beta-\varepsilon}(\mathbb{Q}^\delta)} &\leq \|\tilde{u}^\delta_0\|_{V^1_{-\beta-\varepsilon}(\Omega)} + |a(\delta)|\,\|\tilde{\zeta}\|_{V^1_{-\beta-\varepsilon}(\Omega)} \\
&\leq \|\tilde{u}^\delta_0\|_{V^1_{-\beta-\varepsilon}(\Omega)} + c\,(\,|c^\delta_0| + |C^\delta_0|\,) \\
&\leq c\,(\|f^\delta_{\text{ext}}\|_{\mathring{V}^1_{\beta+\varepsilon}(\Omega)^*} + \|F^\delta_{\text{in}}\|_{\mathring{V}^1_{-\beta+\varepsilon}(\Xi)^*}).
\end{aligned} \tag{38}$$

Now, we bound the term of (37) associated with $U_{\text{in}}^\delta - M^\delta$. On $\mathbb{Q}^{1/\delta}$, thanks to (27), we find $U_{\text{in}}^\delta - M^\delta = \tilde{U}_0^\delta + A(\delta)\tilde{Z} \in \mathring{\mathcal{V}}_\gamma^1(\Xi)$ for all $\gamma \in (0,2)$ ($\tilde{U}_0^\delta$ and $\tilde{Z}$ are respectively defined in (20) and (22)). Therefore, we can write

$$
\begin{aligned}
\|U_{\text{in}}^\delta - M^\delta\|_{\mathcal{V}_{\beta-\varepsilon}^1(\mathbb{Q}^{1/\delta})} &\leq \|\tilde{U}_0^\delta\|_{\mathcal{V}_{\beta-\varepsilon}^1(\Xi)} + |A(\delta)|\,\|\tilde{Z}\|_{\mathcal{V}_{\beta-\varepsilon}^1(\Xi)} \\
&\leq \|\tilde{U}_0^\delta\|_{\mathcal{V}_{\beta-\varepsilon}^1(\Xi)} + c\,(\,|c_0^\delta| + |C_0^\delta|\,) \\
&\leq c\,(\|f_{\text{ext}}^\delta\|_{\mathring{V}_{\beta+\varepsilon}^1(\Omega)^*} + \|F_{\text{in}}^\delta\|_{\mathring{\mathcal{V}}_{-\beta+\varepsilon}^1(\Xi)^*}).
\end{aligned}
\tag{39}
$$

To sum up, plugging (38) and (39) into (37), we find

$$
\|\mathrm{K}^\delta f\|_{\mathring{V}_{\beta,\delta}^1(\Omega)^*} \leq c\,(\delta^\varepsilon + \delta^{\beta-\varepsilon})\,(\|f_{\text{ext}}^\delta\|_{\mathring{V}_{\beta+\varepsilon}^1(\Omega)^*} + \|F_{\text{in}}^\delta\|_{\mathring{\mathcal{V}}_{-\beta+\varepsilon}^1(\Xi)^*}).
\tag{40}
$$

To conclude the estimate, there only remains to examine separately each of the terms appearing on the right hand side above. First, we need a bound for $f_{\text{ext}}^\delta$. Pick an arbitrary $v \in \mathring{V}_{\beta+\varepsilon}^1(\Omega)$ and note that $v\chi_{\sqrt{\delta}} \in \mathrm{H}_0^1(\Omega)$. Observe also that $\|v\chi_{\sqrt{\delta}}\|_{\mathrm{V}_{\beta,\delta}^1(\Omega)} \leq c\,\delta^{-\varepsilon/2}\|v\|_{\mathrm{V}_{\beta+\varepsilon}^1(\Omega)}$. According to (11), we have

$$
\begin{aligned}
|\langle f_{\text{ext}}^\delta, v\rangle_\Omega| = |\langle f, v\,\chi_{\sqrt{\delta}}\rangle_\Omega| &\leq \|f\|_{\mathring{V}_{\beta,\delta}^1(\Omega)^*}\,\|v\chi_{\sqrt{\delta}}\|_{\mathrm{V}_{\beta,\delta}^1(\Omega)} \\
&\leq c\,\delta^{-\varepsilon/2}\|f\|_{\mathring{V}_{\beta,\delta}^1(\Omega)^*}\,\|v\|_{\mathrm{V}_{\beta+\varepsilon}^1(\Omega)}.
\end{aligned}
$$

Dividing by $\|v\|_{\mathrm{V}_{\beta+\varepsilon}^1(\Omega)}$, and taking the supremum over all such $v$, we obtain

$$
\|f_{\text{ext}}^\delta\|_{\mathring{V}_{\beta+\varepsilon}^1(\Omega)^*} \leq c\,\delta^{-\varepsilon/2}\|f\|_{\mathring{V}_{\beta,\delta}^1(\Omega)^*}.
\tag{41}
$$

To derive a bound for $F_{\text{in}}^\delta$, pick an arbitrary $v \in \mathring{\mathcal{V}}_{-\beta+\varepsilon}^1(\Xi)$. Taking into account the definition of $F_{\text{in}}^\delta$ given by (11), and since $\rho \leq 2/\sqrt{\delta}$ in $\mathrm{supp}(\psi_{1/\sqrt{\delta}})$, we see that

$$
\begin{aligned}
|\langle F_{\text{in}}^\delta, v\rangle_\Xi| &\leq \|f\|_{\mathring{V}_{\beta,\delta}^1(\Omega)^*}\|(\tau_\delta \cdot v)\psi_{\sqrt{\delta}}\|_{\mathrm{V}_{\beta,\delta}^1(\Omega)} \\
&\leq \delta^\beta \|f\|_{\mathring{V}_{\beta,\delta}^1(\Omega)^*}\|v\,\psi_{1/\sqrt{\delta}}\|_{\mathcal{V}_\beta^1(\Xi)} \\
&\leq c\,\delta^{\varepsilon/2}\,\|f\|_{\mathring{V}_{\beta,\delta}^1(\Omega)^*}\,\|v\|_{\mathcal{V}_{-\beta+\varepsilon}^1(\Xi)}.
\end{aligned}
$$

Dividing by $\|v\|_{\mathcal{V}_{-\beta+\varepsilon}^1(\Xi)}$ and taking the supremum over all such $v$, we can finally write

$$
\|F_{\text{in}}^\delta\|_{\mathring{\mathcal{V}}_{-\beta+\varepsilon}^1(\Xi)^*} \leq c\,\delta^{\varepsilon/2}\,\|f\|_{\mathring{V}_{\beta,\delta}^1(\Omega)^*}\,.
\tag{42}
$$

At this stage, it is important to notice that in (41) and (42), the exponents obtained for $\delta$ involve $\varepsilon/2$ and not $\varepsilon$. We conclude the proof of (34) by plugging (41) and (42) into (40). To prove that $\mathrm{Id} + \mathrm{K}^\delta$ is invertible for $\delta$ small enough, it suffices to observe that $\|\mathrm{Id}, \mathring{V}_{\beta,\delta}^1(\Omega)^* \to \mathring{V}_{\beta,\delta}^1(\Omega)^*\| = 1$. □

Hence $\hat{\mathrm{R}}^\delta \cdot (\mathrm{Id} + \mathrm{K}^\delta)^{-1}$ is the resolvent of the operator $\mathrm{A}^\delta$. This is an interesting result because it shows that, as long as Assumptions 1 and 2 are satisfied, then Problem (1) is systematically uniquely solvable for $\delta \in \mathrm{I}(\alpha)$ small enough.

## 5.2 Norm estimate of the resolvent

We have established that $(\mathrm{A}^\delta)^{-1} = \hat{\mathrm{R}}^\delta \cdot (\mathrm{Id} + \mathrm{K}^\delta)^{-1}$. Besides Proposition 5.1 allows to obtain an estimate for $(\mathrm{Id} + \mathrm{K}^\delta)^{-1}$. As a consequence, to obtain an estimate for $(\mathrm{A}^\delta)^{-1}$ it suffices to derive a bound for $\hat{\mathrm{R}}^\delta$.

**Proposition 5.2.**
*For any continuous operator $T : \mathrm{H}^{-1}(\Omega) \to \mathrm{H}_0^1(\Omega)$, consider the following $\delta$-dependent operator norm,*

$$\|T, \mathring{\mathrm{V}}_{\beta,\delta}^1(\Omega)^* \to \mathrm{H}_0^1(\Omega)\| := \sup_{g \in \mathrm{H}^{-1}(\Omega) \setminus \{0\}} \frac{\|Tg\|_{\mathrm{H}_0^1(\Omega)}}{\|g\|_{\mathring{\mathrm{V}}_{\beta,\delta}^1(\Omega)^*}}.$$

*Then for any $\beta \in (0,2)$, there exists a constant $c > 0$ independent of $\delta$ such that*

$$\|\hat{\mathrm{R}}^\delta, \mathring{\mathrm{V}}_{\beta,\delta}^1(\Omega)^* \to \mathrm{H}_0^1(\Omega)\| \leq c |\ln \delta|^{1/2}, \qquad \forall \delta \in (0,1) \cap \mathrm{I}(\alpha). \tag{43}$$

**Proof:** For some fixed $f \in \mathrm{H}^{-1}(\Omega)$, we come back to the expression of $\hat{\mathrm{R}}^\delta f = \hat{u}^\delta$ given by (29). Denote $\mathbb{T}^\delta := \{\mathbf{x}(r,\theta) \in \mathbb{R} \times (0,+\infty) \mid r \geq \delta\}$. It is easy to check that $\|\chi_\delta v\|_{\mathrm{H}_0^1(\Omega)} \leq c \|v\|_{\mathrm{V}_0^1(\mathbb{T}^\delta)}$, $\forall v \in \mathrm{H}_0^1(\Omega)$, $\forall \delta \in (0,1)$. Consequently, we have

$$\|\hat{u}^\delta\|_{\mathrm{H}_0^1(\Omega)} \leq c \|u_{\mathrm{ext}}^\delta\|_{\mathrm{V}_0^1(\mathbb{T}^\delta)} + c \|m^\delta\|_{\mathrm{V}_0^1(\mathbb{T}^\delta)} + \|\tau_\delta \cdot (\psi_{1/\delta} U_{\mathrm{in}}^\delta)\|_{\mathrm{H}_0^1(\Omega)}. \tag{44}$$

Use a straightforward calculus, take into account (28), apply (8) (resp. (10)) to $c_0^\delta$ (resp. $C_0^\delta$), and use (41)–(42) with $\varepsilon = 0$, to obtain

$$\begin{aligned}
\|m^\delta\|_{\mathrm{V}_0^1(\mathbb{T}^\delta)} &\leq c |\ln \delta|^{1/2} (|c_0^\delta| + |C_0^\delta|) &\leq& c |\ln \delta|^{1/2} (\|f_{\mathrm{ext}}^\delta\|_{\mathring{\mathrm{V}}_\beta^1(\Omega)^*} + \|F_{\mathrm{in}}^\delta\|_{\mathring{\mathcal{V}}_{-\beta}^1(\Xi)^*}) \\
&&\leq& c |\ln \delta|^{1/2} \|f\|_{\mathring{\mathrm{V}}_{\beta,\delta}^1(\Omega)^*}.
\end{aligned} \tag{45}$$

As $\Omega$ is bounded, there holds $\|v\|_{\mathrm{V}_0^1(\mathbb{T}^\delta \cap \Omega)} \leq c \|v\|_{\mathrm{V}_{-\beta}^1(\Omega)}$ for any $v \in \mathring{\mathrm{V}}_{-\beta}^1(\Omega)$. Besides, $\|\zeta\|_{\mathrm{V}_0^1(\mathbb{T}^\delta)} = O(|\ln \delta|^{1/2})$. As a consequence, using (8) and working like in (45), we find

$$\begin{aligned}
\|u_{\mathrm{ext}}^\delta\|_{\mathrm{V}_0^1(\mathbb{T}^\delta)} &\leq& \|u_0^\delta\|_{\mathrm{V}_0^1(\mathbb{T}^\delta)} + |a(\delta)| \|\zeta\|_{\mathrm{V}_0^1(\mathbb{T}^\delta)} \\
&\leq& c |\ln \delta|^{1/2} (|c_0^\delta| + |C_0^\delta| + \|\tilde{u}_0^\delta\|_{\mathrm{V}_{-\beta}^1(\Omega)}) \\
&\leq& c |\ln \delta|^{1/2} (\|f_{\mathrm{ext}}^\delta\|_{\mathring{\mathrm{V}}_\beta^1(\Omega)^*} + \|F_{\mathrm{in}}^\delta\|_{\mathring{\mathcal{V}}_{-\beta}^1(\Xi)^*}) \leq c |\ln \delta|^{1/2} \|f\|_{\mathring{\mathrm{V}}_{\beta,\delta}^1(\Omega)^*}.
\end{aligned} \tag{46}$$

We finally deal with the third term on the right hand side of (44). Note that $\Xi \cap \mathrm{supp}(\psi_{1/\delta}) = \Xi \setminus \mathbb{T}^{2/\delta}$. We assume that $\delta$ is small enough *i.e.* $\delta \leq 2$, so that $\delta \leq 2/(1+\rho)$ in $\Xi \setminus \mathbb{T}^{2/\delta}$. Then, there holds

$$\begin{aligned}
\|\tau_\delta \cdot (\psi_{1/\delta} U_{\mathrm{in}}^\delta)\|_{\mathrm{H}_0^1(\Omega)} &\leq& c \|U_{\mathrm{in}}^\delta\|_{\mathcal{V}_0^1(\Xi \setminus \mathbb{T}^{2/\delta})} \\
&\leq& c \|U_0^\delta\|_{\mathcal{V}_0^1(\Xi \setminus \mathbb{T}^{2/\delta})} + c |A(\delta)| \|Z(\delta)\|_{\mathcal{V}_0^1(\Xi \setminus \mathbb{T}^{2/\delta})} \\
&\leq& c |\ln \delta|^{1/2} (|c_0^\delta| + |C_0^\delta| + \|\tilde{U}_0^\delta\|_{\mathrm{V}_\beta^1(\Omega)}) \\
&\leq& c |\ln \delta|^{1/2} (\|f_{\mathrm{ext}}^\delta\|_{\mathring{\mathrm{V}}_\beta^1(\Omega)^*} + \|F_{\mathrm{in}}^\delta\|_{\mathring{\mathcal{V}}_{-\beta}^1(\Xi)^*}) \\
&\leq& c |\ln \delta|^{1/2} \|f\|_{\mathring{\mathrm{V}}_{\beta,\delta}^1(\Omega)^*}.
\end{aligned} \tag{47}$$

In the calculus above, we used that $\|Z(\delta)\|_{\mathcal{V}_0^1(\Xi \setminus \mathbb{T}^{2/\delta})} = O(|\ln \delta|^{1/2})$, as well as decomposition (20). To conclude, there remains to plug (45), (46), (47) into (44). Since $\hat{u}^\delta = \hat{\mathrm{R}}^\delta f$, and the estimate (44) holds for any $f \in \mathrm{H}^{-1}(\Omega)$, this finishes the proof. $\square$

Proposition 5.2 was the last building block required in order to establish a norm estimate for the resolvent of $\mathrm{A}^\delta$. We state and prove it now.

**Theorem 5.1.**
*Let $\beta \in (0,2)$. Under Assumptions 1 and 2, there exists $\delta_0$ such that Problem (1) is uniquely solvable for all $\delta \in (0,\delta_0] \cap \mathrm{I}(\alpha)$, with $\alpha \in (0,1/2)$. Moreover, there exists a constant $c > 0$ independent of $\delta$ such that*

$$\|u^\delta\|_{\mathrm{H}_0^1(\Omega)} \leq c |\ln \delta|^{1/2} \|\mathrm{div}(\sigma^\delta \nabla u^\delta)\|_{\mathring{\mathrm{V}}_{\beta,\delta}^1(\Omega)^*}, \qquad \forall u^\delta \in \mathrm{H}_0^1(\Omega), \ \forall \delta \in (0,\delta_0] \cap \mathrm{I}(\alpha).$$

**Proof:** The invertibility result has already been established at the end of §5.1. Here, we have to derive an upper bound for $\|(A^\delta)^{-1}, \mathring{V}^1_{\beta,\delta}(\Omega)^* \to H^1_0(\Omega)\|$ where $A^\delta$ has been defined in (31). According to (32), (34) and (43), one has

$$\|(A^\delta)^{-1}, \mathring{V}^1_{\beta,\delta}(\Omega)^* \to H^1_0(\Omega)\|$$
$$\leq \|\hat{R}^\delta, \mathring{V}^1_{\beta,\delta}(\Omega)^* \to H^1_0(\Omega)\| \cdot \|(\text{Id} + K^\delta)^{-1}, \mathring{V}^1_{\beta,\delta}(\Omega)^* \to \mathring{V}^1_{\beta,\delta}(\Omega)^*\| \leq \frac{c\,|\ln\delta|^{1/2}}{1-(\delta^{\varepsilon/2} + \delta^{\beta - \varepsilon/2})}.$$

□

Admittedly Theorem 5.1 does not yield a bound for $(A^\delta)^{-1}$ that is uniform with respect to $\delta$. It predicts that $\|(A^\delta)^{-1}, \mathring{V}^1_{\beta,\delta}(\Omega)^* \to \mathring{V}^1_{\beta,\delta}(\Omega)^*\|$ does not grow faster than $O(|\ln\delta|^{1/2})$. We shall see in the next section that this does not prevent from obtaining an error estimate concerning an asymptotic expansion of the solution $u^\delta$ to (1).

## 6 First order asymptotics

According to Theorem 5.1, we know that for all $f \in H^{-1}(\Omega)$, Problem (1) is uniquely solvable in the usual functional framework $H^1_0(\Omega)$ for $\delta \in I(\alpha)$ small enough. To obtain this result, we have used a well-suited asymptotic expansion at the first order and adapted weighted spaces. Thus, rounding the corner, with an interface meeting the boundary perpendicularly, could be seen as a regularization method for Problem (4) set in the outer limit geometry. In this section, our goal is to show that this regularization process is not satisfactory. Indeed, we are going to prove that the solution $u^\delta$ depends critically on the value of the parameter $\delta$. To do this, we build a more usual asymptotic expansion of $u^\delta$, with terms that do not depend on $\delta$ (except, of course, the jauge functions), at the first order.

We will assume that the source term $f$ in (1) belongs to $\mathring{V}^1_\beta(\Omega)^*$ with $\beta \in (0, 2)$. Since $H^1_0(\Omega) \subset \mathring{V}^1_\beta(\Omega)$, there holds $\mathring{V}^1_\beta(\Omega)^* \subset H^{-1}(\Omega)$. This implies in particular that $\|f\|_{\mathring{V}^1_{\beta,\delta}(\Omega)^*}$ remains bounded uniformly with respect to $\delta$.

### 6.1 Construction of the expansion

To construct an expansion of $u^\delta$ for $\delta \in I(\alpha)$, we work as in Section 4, using the method of matched asymptotics, but this time, we do not decompose the source term $f$ into outer and inner contributions. This leads us to introduce the function $\check{u}^\delta$ such that

$$\check{u}^\delta(r,\theta) = u_{\text{ext}}(r,\theta)\chi_\delta(r) + U_{\text{in}}(r/\delta,\theta)\psi(r) - m(r,\theta)\chi_\delta(r)\psi(r)$$

$$\text{with } \begin{vmatrix} u_{\text{ext}}(r,\theta) & := & u^0(r,\theta) + a(\delta)\zeta(r,\theta), \\ U_{\text{in}}(\rho,\theta) & := & A(\delta)Z(\rho,\theta), \\ m(r,\theta) & := & A(\delta)(r/\delta)^{i\mu}\phi(\theta) + a(\delta)r^{-i\mu}\phi(\theta)\,. \end{vmatrix} \quad (48)$$

Again, we emphasize that the definition that we provide for $\check{u}^\delta$ differs from the one of $\hat{u}^\delta$ given in (29) because $\check{u}^\delta$ is built from terms that do not depend on $\delta$. In (48), the functions $\zeta(r,\theta), Z(\rho,\theta), \phi(\theta)$ are respectively set in (15)-(16), (19)-(22) and (6). We choose $u^0$ such that $u^0 := (A^{\text{out}}_\beta)^{-1} f \in V^{\text{out}}_\beta(\Omega)$. As a consequence, by the very definition of $V^{\text{out}}_\beta(\Omega)$, there exists a constant $c^0 \in \mathbb{C}$ such that $u^0(r,\theta) - c^0\chi(r)r^{i\mu}\phi(\theta) \in \mathring{V}^1_{-\beta}(\Omega)$. Finally, proceeding as in §4.3, we find the that the jauge functions must verify

$$a(\delta) := \frac{c^0\,C_z\,\delta^{2i\mu}}{1-(\delta/\delta_\star)^{2i\mu}} \qquad \text{and} \qquad A(\delta) := \frac{c^0\,\delta^{i\mu}}{1-(\delta/\delta_\star)^{2i\mu}}\,, \quad (49)$$

where the constants $C_z$ and $\delta_\star$ are defined in (22) and in §4.3. In particular, (49) makes sense thanks to the assumption $\delta \in I(\alpha)$.

### 6.2 Error estimate

Now, computing the error $\|u^\delta - \check{u}^\delta\|_{H^1_0(\Omega)}$, we prove that the function $\check{u}^\delta$ defined by (48) is indeed a sharp approximation of $u^\delta$. First of all, plugging the expression of $u^\delta - \check{u}^\delta$ in the estimate of Theorem 5.1, we obtain, for $\varepsilon \in (0, \beta)$ and $\delta_0$ small enough,

$$\|u^\delta - \check{u}^\delta\|_{H^1_0(\Omega)} \le c\,|\ln \delta|^{1/2}\,\|\mathrm{div}(\sigma^\delta \nabla(u^\delta - \check{u}^\delta))\|_{\mathring{V}^1_{\varepsilon,\delta}(\Omega)^*}, \qquad \forall \delta \in (0, \delta_0] \cap I(\alpha). \tag{50}$$

Let us bound the right hand side of the previous equation. For any $v \in H^1_0(\Omega)$, following §5.1, we can write

$$\begin{aligned}
(\sigma^\delta \nabla \check{u}^\delta, \nabla v)_\Omega &= (\sigma^0 \nabla(u_{\text{ext}} \chi_\delta), \nabla v)_\Omega + (\sigma^\delta \nabla(\tau_\delta \cdot (U_{\text{in}} \psi_{1/\delta})), \nabla v)_\Omega - (\sigma^\delta \nabla(m\,\chi_\delta \psi), \nabla v)_\Omega \\
&= \langle f, \chi_\delta\, v\rangle_\Omega + (\sigma^0 \nabla \chi_\delta, \overline{u}_{\text{ext}} \nabla v - v \nabla \overline{u}_{\text{ext}})_\Omega \\
&\quad + (\sigma^\infty \nabla \psi_{1/\delta}, \overline{U}_{\text{in}} \nabla(\tau_{1/\delta} \cdot v) - (\tau_{1/\delta} \cdot v) \nabla \overline{U}_{\text{in}})_\Xi \\
&\quad - (\sigma^\delta \nabla(\chi_\delta\,\psi), \overline{m}\nabla v - v\nabla \overline{m})_\Omega.
\end{aligned}$$

Since $-\chi_\delta\,\psi = 1 - \chi_\delta - \psi$, denoting, $M := \tau_{1/\delta} \cdot m$, we obtain

$$\begin{aligned}
(\sigma^\delta \nabla(u^\delta - \check{u}^\delta), \nabla v)_\Omega =\ & \langle f, \psi_\delta\, v\rangle_\Omega \\
& + (\sigma^0 \nabla \chi_\delta, (\overline{u}_{\text{ext}} - \overline{m})\nabla v - v\nabla(\overline{u}_{\text{ext}} - \overline{m}))_\Omega \\
& + (\sigma^\infty \nabla \psi_{1/\delta}, (\overline{U}_{\text{in}} - \overline{M})\nabla(\tau_{1/\delta} \cdot v) - (\tau_{1/\delta} \cdot v)\nabla(\overline{U}_{\text{in}} - \overline{M}))_\Xi.
\end{aligned} \tag{51}$$

Let us study each of the three terms of the right hand side of (51). For the first one, there holds,

$$\begin{aligned}
|\langle f, \psi_\delta\, v\rangle_\Omega| \le c\,\|f\|_{\mathring{V}^1_\beta(\Omega)^*}\|\psi_\delta\, v\|_{V^1_\beta(\Omega)} &\le\ c\,\delta^{\beta-\varepsilon}\,\|f\|_{\mathring{V}^1_\beta(\Omega)^*}\|\psi_\delta\,v\|_{V^1_\varepsilon(\Omega)} \\
&\le\ c\,\delta^{\beta-\varepsilon}\,\|f\|_{\mathring{V}^1_\beta(\Omega)^*}\|v\|_{V^1_{\varepsilon,\delta}(\Omega)}.
\end{aligned} \tag{52}$$

For the second and third terms of (51), working as in (36), (38) and (39), one finds successively

$$\begin{aligned}
&\left|\left(\sigma^0\nabla\chi_\delta, (\overline{u}_{\text{ext}} - \overline{m})\nabla v - v\nabla(\overline{u}_{\text{ext}} - \overline{m})\right)_\Omega\right| \\
&\le c\,\|u_{\text{ext}} - m\|_{V^1_{-\beta}(\mathbb{Q}^\delta)}\|v\|_{V^1_\beta(\mathbb{Q}^\delta)} \le c\,\delta^{\beta-\varepsilon}\,\|u_{\text{ext}} - m\|_{V^1_{-\beta}(\mathbb{Q}^\delta)}\|v\|_{V^1_{\varepsilon,\delta}(\Omega)} \\
&\le c\,\delta^{\beta-\varepsilon}\,(|c^0| + \|\tilde{u}^0\|_{V^1_{-\beta}(\Omega)})\|v\|_{V^1_{\varepsilon,\delta}(\Omega)} \le c\,\delta^{\beta-\varepsilon}\,\|f\|_{\mathring{V}^1_\beta(\Omega)^*}\|v\|_{V^1_{\varepsilon,\delta}(\Omega)}
\end{aligned} \tag{53}$$

and

$$\begin{aligned}
&\left|\left(\sigma^\infty\nabla\psi_{1/\delta}, (\overline{U}_{\text{in}} - \overline{M})\nabla(\tau_{1/\delta}\cdot v) - (\tau_{1/\delta}\cdot v)\nabla(\overline{U}_{\text{in}} - \overline{M})\right)_\Xi\right| \\
&\le c\,\|U_{\text{in}} - M\|_{\mathcal{V}^1_\beta(\mathbb{Q}^{1/\delta})}\|\tau_{1/\delta}\cdot v\|_{\mathcal{V}^1_{-\beta}(\mathbb{Q}^{1/\delta})} \\
&\le c\,\delta^\beta\,\|U_{\text{in}} - M\|_{\mathcal{V}^1_\beta(\mathbb{Q}^{1/\delta})}\|v\|_{V^1_{\varepsilon,\delta}(\Omega)} \\
&\le c\,\delta^\beta\,|c^0|\|v\|_{V^1_{\varepsilon,\delta}(\Omega)} \le c\,\delta^\beta\,\|f\|_{\mathring{V}^1_\beta(\Omega)^*}\|v\|_{V^1_{\varepsilon,\delta}(\Omega)}.
\end{aligned} \tag{54}$$

Plugging (52), (53) and (54) in (51), we can divide the resulting inequality by $\|v\|_{V^1_{\varepsilon,\delta}(\Omega)}$ and then take the supremum over all $v \in H^1_0(\Omega)$. This leads to the existence of some constant $c > 0$ independent of $\delta$ such that $\|\mathrm{div}(\sigma^\delta\nabla(u^\delta - \check{u}^\delta))\|_{\mathring{V}^1_{\varepsilon,\delta}(\Omega)^*} \le c\,\delta^{\beta-\varepsilon}\|f\|_{\mathring{V}^1_\beta(\Omega)^*}$. Going back to (50), the previous calculus yields the error estimate $\|u^\delta - \check{u}^\delta\|_{H^1_0(\Omega)} \le c\,|\ln \delta|^{1/2}\,\delta^{\beta-\varepsilon}\|f\|_{\mathring{V}^1_\beta(\Omega)^*}$. Since this estimate holds for all $\varepsilon \in (0, \beta)$, the term $|\ln \delta|^{1/2}$ can be removed.

**Theorem 6.1.**
*Let $\beta \in (0,2)$ and $f \in \mathring{V}^1_\beta(\Omega)^*$. Under Assumptions 1 and 2, there exists $\delta_0$ such that Problem (1) is uniquely solvable for all $\delta \in (0,\delta_0] \cap \mathrm{I}(\alpha)$, with $\alpha \in (0,1/2)$. Moreover, the function $\check{u}^\delta \in \mathrm{H}^1_0(\Omega)$ defined in (48) verifies, for all $\varepsilon$ in $(0,\beta)$,*

$$\|u^\delta - \check{u}^\delta\|_{\mathrm{H}^1_0(\Omega)} \leq c\,\delta^{\beta-\varepsilon}\|f\|_{\mathring{V}^1_\beta(\Omega)^*}, \qquad \forall \delta \in (0,\delta_0] \cap \mathrm{I}(\alpha), \tag{55}$$

*where $c > 0$ is a constant independent of $\delta$ and $f$.*

Thanks to Theorem 6.1, we can prove now the instability phenomenon as $\delta$ tends to zero, the main result of the paper.

### 6.3 Conclusion

Let us study the H$^1$-norm of $u^\delta$ when $\delta$ tends to zero. Denote $\mathbb{S}(\sqrt{\delta}) := \{\mathbf{x}(r,\theta) \in \mathbb{R} \times (0,+\infty) \mid \sqrt{\delta} \leq r \leq 1\}$. There holds $\|\check{u}^\delta\|_{\mathrm{H}^1_0(\Omega)} \geq \|\check{u}^\delta\|_{\mathrm{H}^1_0(\mathbb{S}(\sqrt{\delta}))}$. Moreover, according to (48), we know that for $\mathbf{x} \in \mathbb{S}(\sqrt{\delta})$, the function $\check{u}^\delta$ admits the decomposition $\check{u}^\delta = m(r,\theta) + \tilde{\tilde{u}}^\delta$, where $\tilde{\tilde{u}}^\delta$ is a term which remains bounded for the H$^1$-norm when $\delta$ goes to zero. By a direct computation, one can check that we have $\|m\|_{\mathrm{H}^1_0(\mathbb{S}(\sqrt{\delta}))} \geq c|\ln\delta|^{1/2}$ as soon as $c^0 \neq 0$, where $c^0 \in \mathbb{C}$ is the constant such that $u^0(r,\theta) - c^0\chi(r)r^{i\mu}\phi(\theta) \in \mathring{V}^1_{-\beta}(\Omega)$. Thus, if the source term $f$ is such that the solution $u^0$ is strongly singular ($c^0 \neq 0$), then $\|\check{u}^\delta\|_{\mathrm{H}^1_0(\Omega)} \to +\infty$ when $\delta \to 0$. By virtue of (55), we deduce that $\|u^\delta\|_{\mathrm{H}^1_0(\Omega)} \to +\infty$ when $\delta \in \mathrm{I}(\alpha)$ tends to zero. This is not completely surprising since we know that the limit problem for $\delta = 0$ is not well-posed in $\mathrm{H}^1_0(\Omega)$. But maybe the sequence $(u^\delta)$ converges to some limit for some weaker norm, for example, for the L$^2$-norm? Thanks to our asymptotic expansion, we can prove that the answer to this question is no. Indeed, using (48), one notices that on $\Omega \setminus \mathrm{D}(O,2)$, there holds $\check{u}^\delta = u^0 + a(\delta)\zeta$. Since, the function $\delta \mapsto a(\delta)$ does not converge, we deduce that $(\check{u}^\delta)$ does not converge for the L$^2$-norm when $\delta \to 0$. Again, according to (55), it follows that $(u^\delta)$ does not converge for the L$^2$-norm[2]. The reflexion of the strongly oscillating singularity on the rounded corner modifies the solution $u^\delta$ everywhere in $\Omega$, and not only in a neighbourhood of $O$. As a consequence, $u^\delta$ critically depends on the value of $\delta$, it is not stable.

## 7 Numerical illustration

As an illustration of the results obtained in the previous sections, let us conclude by studying on a canonical geometry the question of the stability of Problem (1) with respect to the parameter $\delta$. The geometry will be chosen so that we can separate variables and thus, proceed to explicit computations. Again, we point out the following result. If the limit problem, *i.e.* Problem (1) with $\delta = 0$, is well-posed in $\mathrm{H}^1_0(\Omega^0)$, that boils down to say, for the chosen configuration, that $\kappa_\sigma = \sigma_-/\sigma_+ \in (-\infty,-1) \cup (-1/3,0)$, then the solution $u^\delta$ is well-defined for $\delta$ small enough and converges to $u^0$ for the H$^1$-norm. Conversely, if the limit problem is ill-posed in $\mathrm{H}^1_0(\Omega^0)$ (not of Fredholm type), then Problem (1) critically depends on $\delta$, even for small $\delta$. The framework (see Fig. 6) will be slightly different from the one introduced in §2 because $\overline{\Omega^\delta_+} \cup \overline{\Omega^\delta_-}$ will not be a fixed domain. However, the analysis we provided all along this paper could be extended without difficulty to the geometry studied in this paragraph and we would obtain analogous results. Finally, we present numerical simulations, using a standard finite element approximation, to illustrate the difference in behaviour of the sequence $(u^\delta)$, depending on whether or not the contrast belongs to the critical interval $(-1,-1/3)$.

---
[2]In other words, $(u^\delta)$ and $(\check{u}^\delta)$ are sequences which do not converge, even for the L$^2$-norm, but which are such that, according to (55), $(u^\delta - \check{u}^\delta)$ converges in $\mathrm{H}^1_0(\Omega)$ when $\delta \to 0$.

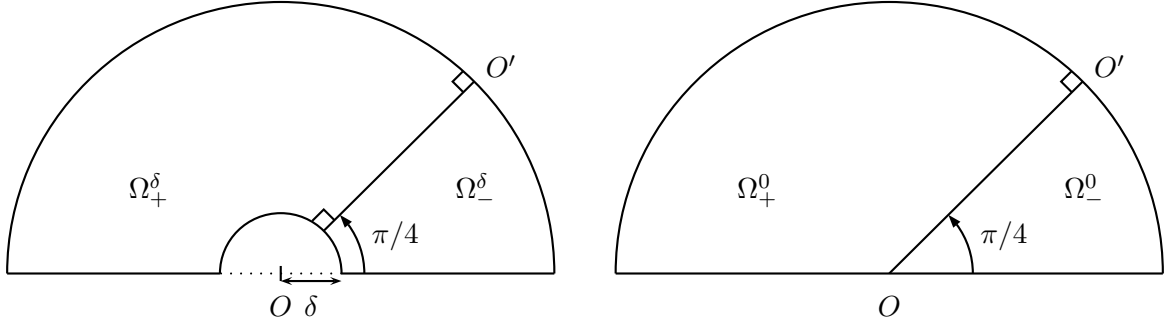

Fig. 6: Domains $\Omega^\delta$ and $\Omega^0$.

Let us first describe the geometry. Consider $\delta \in (0,1)$ and define (see Fig. 6)

$$\begin{array}{rcl}
\Omega^\delta_+ & := & \{\,(r\cos\theta, r\sin\theta)\mid \delta < r < 1,\ \pi/4 < \theta < \pi\,\}; \\
\Omega^\delta_- & := & \{\,(r\cos\theta, r\sin\theta)\mid \delta < r < 1,\ 0 < \theta < \pi/4\,\}; \\
\Omega^\delta & := & \{\,(r\cos\theta, r\sin\theta)\mid \delta < r < 1,\ 0 < \theta < \pi\,\}.
\end{array}$$

Again, introduce the function $\sigma^\delta : \Omega \to \mathbb{R}$ such that $\sigma^\delta = \sigma_\pm$ in $\Omega^\delta_\pm$, where $\sigma_+ > 0$ and $\sigma_- < 0$ are constants. We define the continuous linear operator $A^\delta : H^1_0(\Omega^\delta) \to H^{-1}(\Omega^\delta)$ such that

$$\langle A^\delta u, v\rangle_{\Omega^\delta} = (\sigma^\delta \nabla u, \nabla v)_{\Omega^\delta}, \quad \forall u,v \in H^1_0(\Omega^\delta).$$

Since the interface $\Sigma^\delta := \overline{\Omega^\delta_+} \cap \overline{\Omega^\delta_-}$ is smooth and meets $\partial \Omega^\delta$ orthogonally, we can establish, using the T-coercivity technique (see [10, 4]), that for $\kappa_\sigma = \sigma_-/\sigma_+ \in \mathbb{R}^*\setminus\{-1\}$, for all $\delta \in (0,1)$, $A^\delta$ is a Fredholm operator of index 0. Let us study the question of the injectivity of $A^\delta$.

**T-coercivity approach.** In the sequel, if $v$ is a function defined on $\Omega^\delta$, we will denote $v_+ := v|_{\Omega^\delta_+}$ and $v_- := v|_{\Omega^\delta_-}$. Introduce the operators $R_+$, $R_-$, $T_+$ and $T_-$ such that, for all $u \in H^1_0(\Omega^\delta)$,

$$(R_+ u_+)(r,\theta) = u_+(r, \pi - 3\theta) \quad ; \quad (R_- u_-)(r,\theta) = \begin{cases} u_-(r, \pi/2 - \theta) & \text{if } \theta \le \pi/2 \\ 0 & \text{else} \end{cases};$$

$$T_+ u = \begin{cases} u_+ & \text{on } \Omega^\delta_+ \\ -u_- + 2R_+ u_+ & \text{on } \Omega^\delta_- \end{cases} \quad ; \quad T_- u = \begin{cases} u_+ - 2R_- u_- & \text{on } \Omega^\delta_+ \\ -u_- & \text{on } \Omega^\delta_-. \end{cases}$$

One has $u_+ = -u_- + 2R_+ u_+$ on $\Sigma^\delta$, so $T_+$ is valued in $H^1_0(\Omega^\delta)$. Noticing that $T_+ \cdot T_+ = \text{Id}$, we deduce that $T_+$ is an isomorphism of $H^1_0(\Omega^\delta)$. The same result holds true for $T_-$. For all $u \in H^1_0(\Omega^\delta)$, with the help of Young's inequality, we can write, for all $\eta > 0$,

$$\begin{array}{rcl}
|\langle A^\delta u, T_+ u\rangle_{\Omega^\delta}| & = & |(\sigma_+ \nabla u_+, \nabla u_+)_{\Omega^\delta_+} + (|\sigma_-|\nabla u_-, \nabla u_-)_{\Omega^\delta_-} - 2(|\sigma_-|\nabla u_-, \nabla(R_+ u_+))_{\Omega^\delta_-}| \\
& \ge & (\sigma_+ \nabla u_+, \nabla u_+)_{\Omega^\delta_+} + (|\sigma_-|\nabla u_-, \nabla u_-)_{\Omega^\delta_-} \\
& & -\eta\,(|\sigma_-|\nabla u_-, \nabla u_-)_{\Omega^\delta_-} - 1/\eta\,(|\sigma_-|\nabla(R_+ u_+), \nabla(R_+ u_+))_{\Omega^\delta_-} \\
& \ge & ((\sigma_+ - \|R_+\|^2\,|\sigma_-|/\eta)\,\nabla u_+, \nabla u_+)_{\Omega^\delta_+} + (|\sigma_-|\,(1-\eta)\,\nabla u_-, \nabla u_-)_{\Omega^\delta_-}.
\end{array}$$

Therefore, if $\sigma_+/|\sigma_-| > \|R_+\|^2 = 3$ ($\Leftrightarrow \kappa_\sigma > -1/3$), there exists $c > 0$ independent of $\delta$ such that

$$c\,\|u\|^2_{H^1_0(\Omega^\delta)} \le |\langle A^\delta u, T_+ u\rangle_{\Omega^\delta}|, \quad \forall u \in H^1_0(\Omega^\delta).$$

In the case where $\kappa_\sigma < -1$, we can proceed in the same way, using $T_-$ instead of $T_+$ (notice that $\|R_-\|^2 = 1$). The previous discussion leads to a result of uniform stability under some restrictive conditions on $\kappa_\sigma$.

**Proposition 7.1.**
*Assume that $\kappa_\sigma \in (-\infty, -1) \cup (-1/3, 0)$. Then, for all $\delta \in (0,1)$, the operator $A^\delta$ is an isomorphism. Moreover, there exists $c > 0$ independent of $\delta$ such that*

$$\|u^\delta\|_{H_0^1(\Omega^\delta)} \leq c\, \|A^\delta u^\delta\|_{H^{-1}(\Omega^\delta)}, \qquad \forall u^\delta \in H_0^1(\Omega^\delta).$$

This proposition provides a norm estimate of the resolvent of $A^\delta$ with a constant that does not blow up as $\delta$ tends to 0. But this result holds only for a contrast lying outside the critical interval. For $\kappa_\sigma \in (-1, -1/3)$, the method of the T-coercivity no longer works.

**Direct approach.** Let us go back to the question of the injectivity of $A^\delta$. Instead of using the T-coercivity approach, let us compute explicitly the elements of $\operatorname{Ker} A^\delta$. Consider $u^\delta$ a function of $H_0^1(\Omega^\delta)$ such that $A^\delta u^\delta = 0$. The pair $(u_+^\delta, u_-^\delta)$ satisfies

$$\left| \begin{array}{ll} \Delta u_+^\delta = 0 & \text{in } \Omega_+^\delta \\ \Delta u_-^\delta = 0 & \text{in } \Omega_-^\delta \\ u_+^\delta - u_-^\delta = 0 & \text{on } \Sigma^\delta \\ \sigma_+ \partial_\theta u_+^\delta - \sigma_- \partial_\theta u_-^\delta = 0 & \text{on } \Sigma^\delta. \end{array} \right. \tag{56}$$

With the change of variables $(t, \theta) := (\ln r, \theta)$, Problem (56) is changed in a problem set in the truncated strip $\mathcal{S}^\delta := (\ln \delta, 0) \times (0, \pi)$. In this geometry, separation of variables is a natural approach. Using this trick, it is easy to prove that the family $\{r \mapsto \sin(n\pi(\ln r - \ln \delta)/\ln \delta)\}_{n=1}^\infty$ is a basis of $L^2(\delta, 1)$. Decomposing $u_+^\delta$ and $u_-^\delta$ and since $\Delta u_+^\delta = \Delta u_-^\delta = 0$, one finds

$$u_+^\delta(r, \theta) = \sum_{n=1}^{+\infty} u_{+,n}^\delta \, \sinh\!\left(n\pi \frac{(\theta - \pi)}{\ln \delta}\right) \sin\!\left(n\pi \frac{\ln(r/\delta)}{\ln \delta}\right)$$

$$\text{and} \quad u_-^\delta(r, \theta) = \sum_{n=1}^{+\infty} u_{-,n}^\delta \, \sinh\!\left(n\pi \frac{\theta}{\ln \delta}\right) \sin\!\left(n\pi \frac{\ln(r/\delta)}{\ln \delta}\right),$$

where $u_{+,n}^\delta$ and $u_{-,n}^\delta$ are some constants. The transmission conditions lead to the additional relations, for all $n \in \mathbb{N}^*$,

$$\left| \begin{array}{rcl} -u_{+,n}^\delta \sinh\left(3n\pi^2/(4\ln\delta)\right) & = & u_{-,n}^\delta \sinh\left(n\pi^2/(4\ln\delta)\right) \\ u_{+,n}^\delta \sigma_+ \cosh\left(3n\pi^2/(4\ln\delta)\right) & = & u_{-,n}^\delta \sigma_- \cosh\left(n\pi^2/(4\ln\delta)\right). \end{array} \right. \tag{57}$$

For $n \in \mathbb{N}^*$, the system (57) with respect to $(u_{+,n}^\delta, u_{-,n}^\delta)$ has a non trivial solution if and only if its determinant vanishes which yields the equations, with the notation $\nu_n^\delta := n\pi^2/(4\ln\delta)$,

$$\sigma_- \sinh(3n\pi^2/(4\ln\delta))\cosh(n\pi^2/(4\ln\delta)) + \sigma_+ \sinh(n\pi^2/(4\ln\delta))\cosh(3n\pi^2/(4\ln\delta)) = 0$$

$$\Leftrightarrow \quad \sigma_-\left(\sinh(4\nu_n^\delta) + \sinh(2\nu_n^\delta)\right) + \sigma_+\left(\sinh(4\nu_n^\delta) - \sinh(2\nu_n^\delta)\right) = 0$$

$$\Leftrightarrow \quad \cosh(2\nu_n^\delta) = \frac{\sigma_+ - \sigma_-}{2(\sigma_+ + \sigma_-)} = \frac{1 - \kappa_\sigma}{2(1 + \kappa_\sigma)}.$$

• CASE $\kappa_\sigma \in (-\infty, -1) \cup (-1/3, 0)$
For this case, there holds $(1-\kappa_\sigma)/(2(1+\kappa_\sigma)) < 1$. Consequently, for all $n \in \mathbb{N}^*$, the only solution of (57) is the null solution and $A^\delta$ is injective. This is consistent with Proposition 7.1.

• CASE $\kappa_\sigma \in (-1, -1/3)$
In this situation, one has $(1-\kappa_\sigma)/(2(1+\kappa_\sigma)) > 1$. For all $n \in \mathbb{N}^*$, there exists a unique $\delta^n$ such that (57) has a non trivial solution:

$$\delta^n = \exp\left(-\frac{n\pi^2}{2\operatorname{acosh}(\frac{1-\kappa_\sigma}{2(1+\kappa_\sigma)})}\right) \xrightarrow[n\to\infty]{} 0. \tag{58}$$

Hence, $A^\delta : H^1_0(\Omega^\delta) \to H^{-1}(\Omega^\delta)$ is an isomorphism if and only if $\delta \in (0,1) \setminus \cup_{n \in \mathbb{N}^*} \{\delta^n\}$. Thus, we obtain by an other method, for this special geometry, the result of Proposition 5.1. Note that for this particular domain $\Omega^0$, we can prove (see [7, Lemm. 4.1]) that Assumptions 1 and 2 hold.

**Numerical illustration.** Let us check this numerically. For the computations, we use the *FreeFem++*[3] software while we display the results with *Matlab*[4] and *Paraview*[5]. For details concerning the discretization of Problem (1), we refer the reader to [10, 40, 18]. Take $f \in L^2(\Omega)$ whose support does not meet $O$. Here, we choose $f$ such that $f(x,y) = 100$ if $x < 0.5$ and $f(x,y) = 0$ if $x \geq 0.5$. Moreover, we impose $\sigma_+ = 1$. Let us consider $(\mathcal{T}^\delta_h)_h$ a shape regular family of triangulations of $\overline{\Omega^\delta}$, made of triangles. We assume that, for any triangle $\tau$, one has either $\tau \subset \overline{\Omega^\delta_+}$ or $\tau \subset \overline{\Omega^\delta_-}$. Define the family of finite element spaces

$$V^\delta_h := \left\{ v \in H^1_0(\Omega^\delta) \text{ such that } v|_\tau \in \mathbb{P}_1(\tau) \text{ for all } \tau \in \mathcal{T}^\delta_h \right\},$$

where $\mathbb{P}_1(\tau)$ is the space of polynomials of degree at most 1 on the triangle $\tau$. Let us consider the problem

$$\left| \begin{array}{l} \text{Find } u^\delta_h \in V^\delta_h \text{ such that} \\ (\sigma^\delta \nabla u^\delta_h, \nabla v^\delta_h)_{\Omega^\delta} = (f, v^\delta_h)_{\Omega^\delta}, \quad \forall v^\delta_h \in V^\delta_h. \end{array} \right. \tag{59}$$

• OUTSIDE THE CRITICAL INTERVAL
In Fig. 7, we observe the variation of $\|u^\delta_h\|_{H^1_0(\Omega^\delta)}$ with respect to $1-\delta$ for $\kappa_\sigma = \sigma_-/\sigma_+ = -1-10^{-4}$. Fig. 8 represents the solution $u^\delta_h$ for this contrast and for eight values of $\delta$. As it was expected, Problem (59) seems to be stable with respect to $\delta$.

• INSIDE THE CRITICAL INTERVAL
In Fig. 9, we display the variation of $\|u^\delta_h\|_{H^1_0(\Omega^\delta)}$ with respect to $1-\delta$ for $\kappa_\sigma = \sigma_-/\sigma_+ = -1+10^{-4}$. In accordance with the analytical computations resulting in formula (58), we observe peaks which correspond to the values $\delta = \delta^n$ for which $A^\delta$ fails to be injective. More generally the numerical experiments confirm the following important idea: when the contrast lies inside the critical interval, the solution in the rounded geometry, when it is well-defined, critically depends on the parameter $\delta$. This appears clearly in Fig. 10 where we see that the solution is not stable with respect to $\delta$ not only in a neighbourhood of the corner point $O$ but everywhere in $\Omega$ (cf. discussion of §6.3). Concerning the numerical analysis point of view, let us make the following comments. First, notice that for small values of $\delta$, it is very expensive to use a mesh adapted to the geometry. Therefore, the mesh size is chosen small but independent of $\delta$. This explains why peaks do not appear for $\delta \leq 0.02$. Remark also that we work here with a contrast very close to $-1$. This may seem surprising because for $\kappa_\sigma = -1$, the operators $A^\delta$ are not of Fredholm type, due to the presence of singularities all over the interface [4, Th. 6.2]. However, this allows us to obtain several peaks without being obliged to use a very refined mesh. Indeed, in this case, in (58), the coefficient $\pi^2/(2\operatorname{acosh}(\frac{1-\kappa_\sigma}{2(1+\kappa_\sigma)}))$ is of order $-0.5$. From a numerical point of view, it only requires to use meshes which are locally symmetric with respect to the interface to avoid instability phenomena (see [18]).

---

[3] *FreeFem++*, http://www.freefem.org/ff++/.
[4] *Matlab*, http://www.mathworks.se/.
[5] *Paraview*, http://www.paraview.org/.

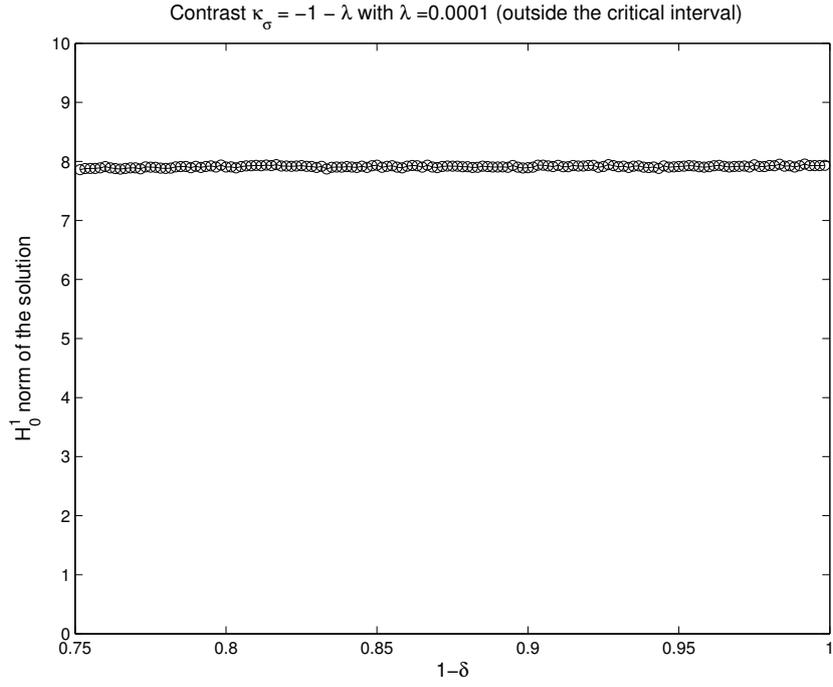

Fig. 7: Variation of $\|u_h^\delta\|_{\mathrm{H}_0^1(\Omega^\delta)}$ with respect to $1 - \delta$.

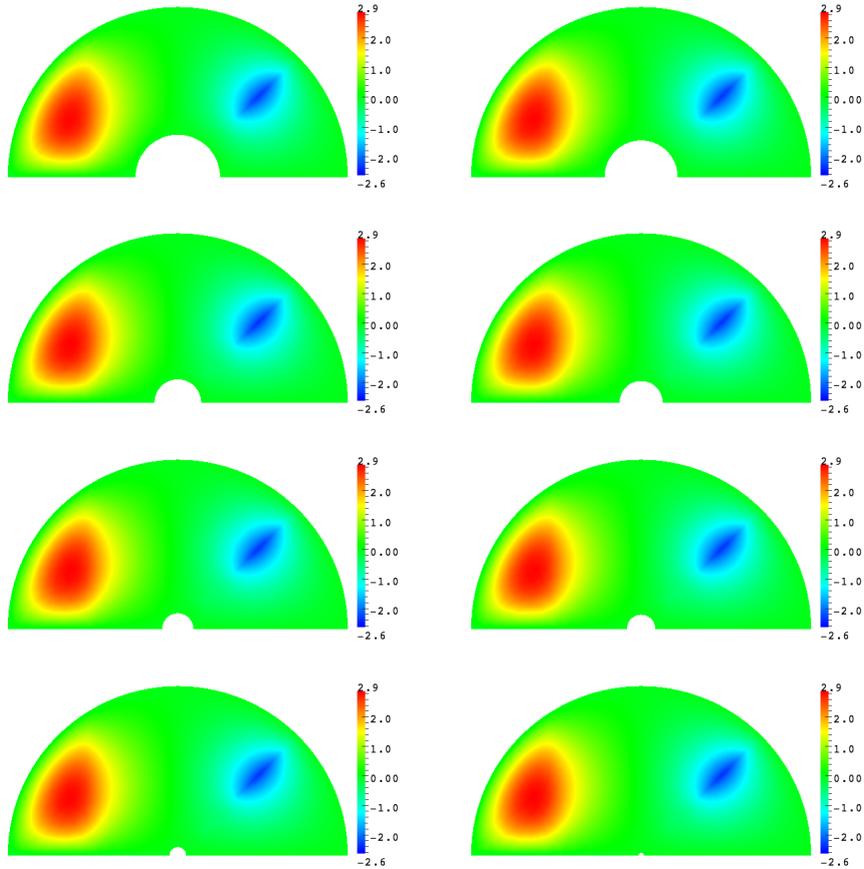

Fig. 8: Solution $u_h^\delta$ for eight values of $\delta$. The contrast $\kappa_\sigma = \sigma_-/\sigma_+$ is chosen equal to $-1 - 10^{-4}$ (outside the critical interval). The sequence $(u_h^\delta)_\delta$ seems to converge when $\delta \to 0$.

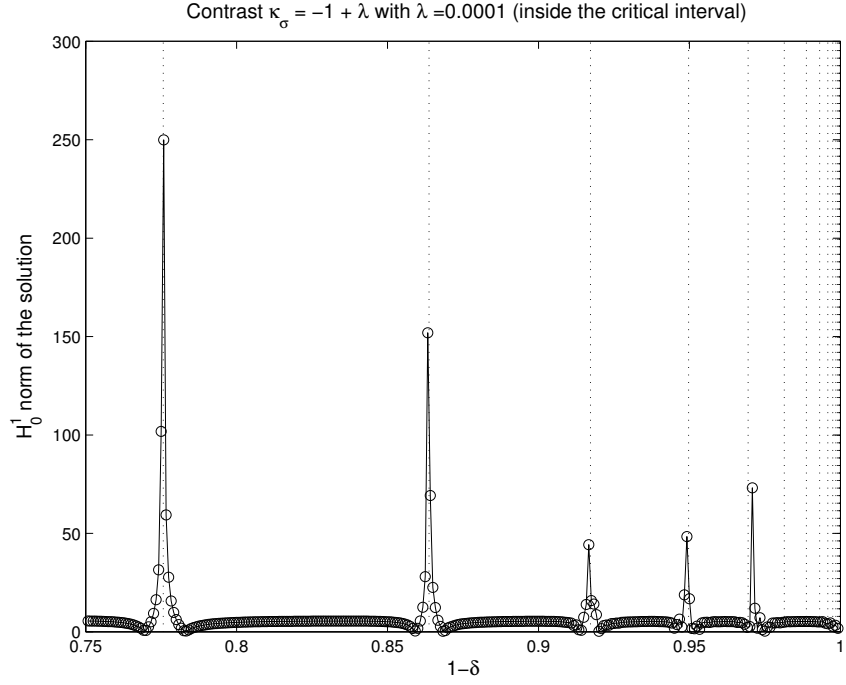

Fig. 9: Variation of $\|u_h^\delta\|_{\mathrm{H}_0^1(\Omega^\delta)}$ with respect to $1-\delta$. The dotted lines correspond to the expected values of $\delta = \delta^n$ for which $A^\delta$ fails to be injective (see the computations resulting in (58)).

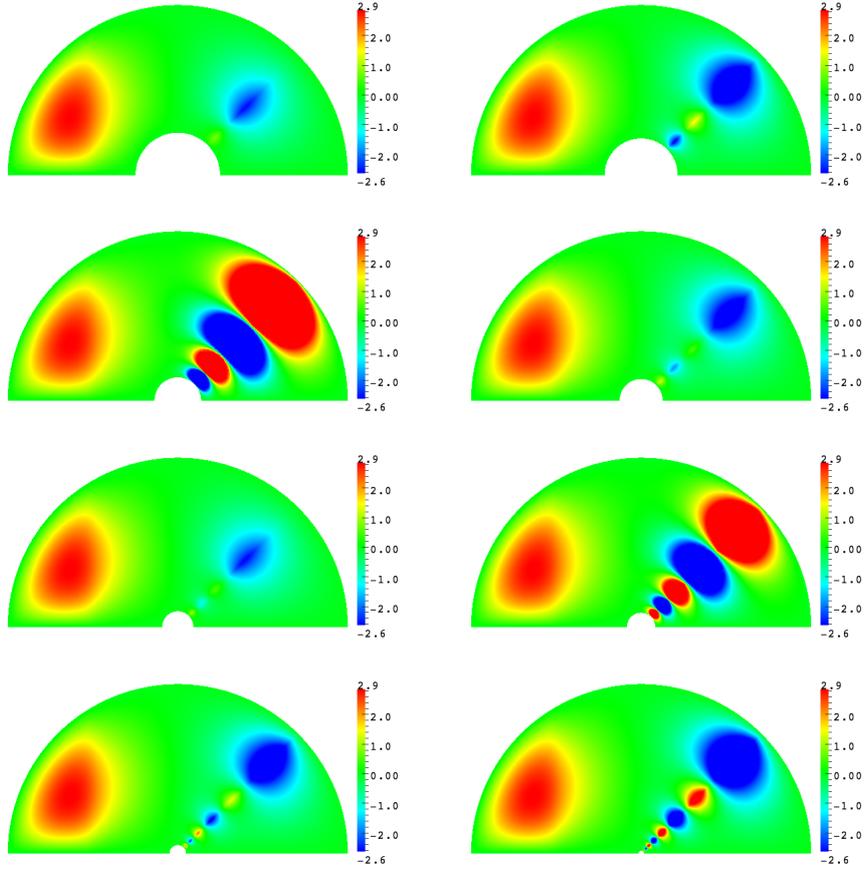

Fig. 10: Solution $u_h^\delta$ for eight values of $\delta$. The contrast $\kappa_\sigma = \sigma_-/\sigma_+$ is chosen equal to $-1 + 10^{-4}$ (inside the critical interval). The sequence $(u_h^\delta)_\delta$ seems not to converge when $\delta \to 0$.

# Acknowledgments

The research of the two first authors is supported by the ANR project METAMATH, grant ANR-11-MONU-016 of the French Agence Nationale de la Recherche. The work of the first author is also supported by the Academy of Finland (decision 140998). The research of the third author is supported by the Russian Foundation for Basic Research, grant No. 12-01-00348.